\documentclass[authoryear, 11pt]{elsarticle}
\usepackage[utf8]{inputenc}
\usepackage{subfiles} 

\usepackage[letterpaper,margin=1in]{geometry}
\usepackage{threeparttable}
\usepackage{booktabs} 
\usepackage{multirow}
\usepackage{graphicx}
\usepackage{makecell} 
\usepackage{amsmath, amssymb}
\usepackage{hyperref}
\usepackage{pdflscape}
\usepackage{tabularx}
\usepackage{array} 
\usepackage{transparent}
\usepackage{longtable}
\newcommand{\tinyscript}{\fontsize{6}{7.2}\selectfont}
\hypersetup{
    pdfauthor={Paola Munoz Briones, Meng-Lin Tsai, Styliani Avraamidou}
}

\usepackage{lineno}
\usepackage{setspace}
\journal{Journal of Cleaner Production}

\makeatletter
\def\ps@pprintTitle{%
 \let\@oddhead\@empty
 \let\@evenhead\@empty
 \def\@oddfoot{}%
 \let\@evenfoot\@oddfoot}
\makeatother

\begin{document}

\begin{frontmatter}

\title{Incorporating circular economy policies into product supply chains using bilevel optimization - A case study on coffee packaging}

\author[a]{Paola Munoz Briones\fnref{equal}}
\author[a]{Meng-Lin Tsai\fnref{equal}}
\author[a]{Styliani Avraamidou\corref{cor1}}

\fntext[equal]{These authors contributed equally to this work.}
\cortext[cor1]{Corresponding author. Email: avraamidou@wisc.edu}

\affiliation[a]{Department of Chemical & Biological Engineering, University of Wisconsin-Madison, Madison, WI, 53706, USA}

\begin{abstract}
Transitioning to a Circular Economy requires policies to drive sustainable practices. This study proposes a bilevel optimization framework to evaluate the combined use of carbon taxes and subsidies in promoting circular supply chains under varying budget levels. A case study of the coffee packaging supply chain with an Extended Producer Responsibility scenario is used to demonstrate this approach. The framework captures the hierarchical interaction between a regional government (upper level), which aims to minimize environmental impacts, and coffee companies (lower level), which seek to minimize costs. Two bilevel optimization problems are formulated based on two environmental objectives: (1) minimization of greenhouse gas (GHG) emissions, and (2) maximization of circularity. The model integrates mixed-integer linear programming (MILP) with life cycle assessment (LCA), techno-economic assessment (TEA) and circularity assessment. Results demonstrate that subsidies effectively drive supply chain shifts toward low-emission and high-circularity configurations, while carbon taxes alone have a more limited impact. Sensitivity analyses highlight the influence of key parameters, such as glass washing distance and loss rates, on policy effectiveness. Overall, the study provides a bilevel optimization framework with quantitative insights to support policy design for sustainable circular supply chains.

\end{abstract}



\end{frontmatter}

\section{Introduction}
 Industrial supply chains involve raw material extraction and processing, which contribute to environmental issues such as resource depletion, greenhouse gas (GHG) emissions, and waste generation. For instance, supply chains producing steel, iron, cement, chemicals, and petrochemicals account for nearly 25\% of global GHG emissions \citep{Teixeira2025}. In the agro-food supply chain, agriculture alone accounts for 70\% of water withdrawals and contributes to water pollution \citep{Sagasta2017}. Consequently, circular economy (CE) emerges as a framework to promote the transition towards achieving environmental, economic, and social sustainability by taking actions to minimize environmental impacts, promote the closed-loop handling of resources, and preserve the value of materials and products for as long as possible \citep{ellenmacarthur2024}. There is a growing need to integrate CE practices within supply chains, including the adoption of innovative manufacturing methods, waste valorization processes, advanced waste management technologies, and the procurement and use of environmentally sustainable materials. 

Supply chains involve multiple stakeholders, including companies and governmental bodies, which often present conflicting objectives; thereby, hindering the alignment of economic and environmental goals. Government-industry collaboration is widely recognized as a key component of CE supply chains \citep{sudusinghe2022supply}. Since the actions of one stakeholder can influence the decisions and outcomes of others, it is essential to consider both perspectives in supply chain modeling \citep{avraamidou2020circular}. For example, a government policy such as a subsidy can shift a firm's technology or process selection, affecting the emissions and profit generated \citep{Jin2014}. This highlights the importance of considering the diverse stakeholder interests in the design and operation of sustainable supply chains. Consequently, evaluating the impact of such policies is essential to ensure their effectiveness.

 Transitioning to environmentally sustainable alternatives often faces economic barriers, as these options may not yet be cost-competitive with conventional practices. Strategies such as policies and regulations play a key role in promoting and accelerating the adoption of sustainable practices across supply chains. Different countries and regions have implemented a broad range of policy measures to accelerate the transition towards CE. Table~\ref{tab:policy_measures} summarizes different policy instruments implemented and provides examples of their implementation. By setting reduction targets, banning specific materials, implementing eco-design regulations, and introducing environmental taxes, governmental agencies try to encourage sustainable practices across supply chain stages—from sourcing to end-of-life. Government subsidies can also influence supply chain management decisions and the adoption of more eco-friendly technologies \citep{li2023examining}. Overall, policies can enable the development of CE supply chains that are both environmentally and economically sustainable. 

\newcolumntype{M}[1]{>{\centering\arraybackslash}m{#1}}

\setlength\extrarowheight{2pt} 
\begin{landscape}
\begin{table}[htbp]
\begin{threeparttable}
\centering
\tinyscript
\caption{Summary of policy measures supporting CE supply chains}
\label{tab:policy_measures}
\begin{tabular}
{M{3.8cm}M{3.8cm}M{2cm}m{9.5cm}m{1cm}}
\hline
\textbf{Policy Measure} & \textbf{Objective} & \textbf{Country / Region} & \textbf{\hspace{4 cm} Policy Description}&\textbf{Year} \\
\hline
\multirow[c]{2}{=}{\centering \\Set reduction targets and regulations} 
& \multirow[c]{2}{=}{\centering \\Reduce overall consumption} 
& Finland 
& Established a cap on the domestic use of natural resources. The overall domestic use of primary raw materials in 2035 will not exceed the 2015 level.\citep{Mikkonen2021Finland} & 2021\\
& & Finland 
& Introduced an energy tax reform, prepared a roadmap and actions to reduce transport emissions by 2030. \citep{Mikkonen2021Finland} &2021 \\
\hline
\multirow[c]{2}{=}{\centering \\Ban specific materials} 
& \multirow[c]{2}{=}{\centering \\Restrict market availability} 
& EU 
& Banned single-use plastics (cotton bud sticks, cutlery, plates, stirrers, straws, cups, food and beverage containers made of expanded polystyrene, and all products made of oxo-degradable plastic).\citep{EU_Directive2019_904} &2021\\
& & Finland 
& Banned the use of coal for power generation after May 2029.\citep{Mikkonen2021Finland} &2019 \\
& & Chile 
& Banned plastic shopping bags, certain single-use plastic items (cutlery, straws, stirrers, expanded polystyrene containers, etc.) in food service establishments. Supermarkets need to offer returnable beverage bottles.\citep{Chile_MMA_EconomiaCircular_Plasticos}
&2022 \\

\hline
\multirow[c]{2}{=}{\centering \\[3pt]Design regulations} 
& \multirow[c]{2}{=}{\centering \\[0.5pt]Require sustainable production design} 
& EU 
& Established resource-efficiency requirements for devices such as refrigerators, dishwashers, electronic displays, washing machines, welding equipment, and data servers/storage products.\citep{CalistoFriant2021} &2019\\
& & EU 
& Required PET beverage bottles to contain at least 25\% recycled plastic starting in 2025, and mandates 30\% recycled plastic by 2030. \citep{EU_Commission2025_SUP} &2019\\
\hline
Provide incentives 
& Encourage eco-design innovation and adoption of sustainable practices
& France 
& Provided tax credits for Green Building and Green Industry (e.g., solar panel manufacturing). \citep{ClarkeVivier2024_FrenchIncentives} &2024\\
& & Canada
&Toronto supports small businesses through the "The Circular Food Innovators Fund" to introduce reusable food service systems to replace single-use and takeaway items.\citep{toronto_circular_food_innovators_fund}
&2024 \\
\hline
\multirow[c]{2}{=}{\centering \\[6pt]Labeling regulations}
& \multirow[c]{2}{=}{\centering Mandate clear product marking on repairing or waste management or environmental impact}
& EU 
& By July 2021, labeling on waste management options and environmental impacts must appear on beverage cups, tobacco, wet wipes, sanitary towels, etc.\citep{EU_Directive2019_904} &2019*\\
& & France & Introduced a Repairability Index for electronic devices, like smartphones, to inform consumers about repair options before purchase.\citep{Microsoft2025_RepairabilityIndexFrance} &2021\\
\hline
\multirow[c]{2}{=}{\centering \\[3pt]Collection regulations} 
& \multirow[c]{2}{=}{\centering \\Set targets for separate waste collection} 
& EU 
& Establishes a target for 77\% separation of plastic bottles by 2025 and 90\% by 2029.\citep{EU_Directive2019_904} &2019\\
& & Germany 
& Mandates households sort waste into designated bins: yellow for lightweight packaging, green/blue for paper/glass, black/grey for residual waste, brown for biowaste.\citep{HandbookGermany_WasteSeparation}&2015 \\
\hline
Awareness campaigns and programs 
& Increase public awareness and engagement 
& Netherlands 
& The ``Waste at School" program focuses on reducing waste generation within schools.\citep{EEA2023_NetherlandsWastePrevention} &2021 \\
\hline
\multirow[c]{2}{=}{\centering \\[4pt] Environmental taxes}
& \multirow[c]{2}{=}{\centering \\[0.5pt] Internalize environmental costs} 
& Germany 
& Introduced a national ETS applying a carbon price to heating and transport fuels. \citep{ICAP_Germany_nETS2025} &2021\\
& & China 
 & Implement tax policies towards energy conservation, environmental protection and resource utilization to guarantee green development through the ``14th Five-Year Plan" \citep{China_14thFiveYearPlan_2021}
& 2021 \\
\hline
\end{tabular}
\begin{tablenotes}
\scriptsize
\item[*] Enforced in 2021.
\end{tablenotes}
\end{threeparttable}
\end{table}

\end{landscape}

\subsection{Modeling the efficacy and effect of policy instruments on supply chains}
Researchers employ modeling, optimization, and game theory to integrate environmental policies into their supply chain studies, each with differing levels of complexity and detail.

Initially, studies primarily focused on either carbon taxes or incentives, with limited investigation into their interaction. For instance, Jin et al \cite{Jin2014} use optimization models and sensitivity to explore the impacts of different carbon emissions policies (carbon emission tax, inflexible cap, and cap-and-trade) on supply chain design and logistics operation, focusing on emissions and costs. Ma et al \cite{ma2021infrastructures} explore different types of economic subsidies to balance between high-cost subsidy solutions and low-cost conventional polluting practices in the dairy supply chain. Despite these efforts, the models separately examine either carbon taxes or subsidies, overlooking the interactive effects of these two policies.

Game theory approaches have been increasingly used to capture the interactions among supply chain stakeholders in modeling. One method to account for multiple decision-makers and their objectives is to use Nash equilibrium to model the interactions \citep{tominac2017game,torres2016design}, where the players take actions simultaneously until an equilibrium is reached. Another approach is the use of evolutionary game theory models to study how the government's and manufacturer's strategies evolve over time in response to policy dynamics. Chen et al \citep{Chen2018} develop an evolutionary game theory model to examine different combinations of carbon taxes and subsidies to find the optimal mechanism. To capture hierarchical, leader-follower relationships involving sequential decision-making, Stackelberg game models and multi-level optimization are used, where decision-makers take their decisions sequentially \citep{fortuny1981representation}.

\subsection{Stackelberg games and bilevel optimization for CE supply chains}

A Stackelberg game, models sequential decision-making between two players: a leader and a follower \citep{von1952theory}. The leader takes a move (i.e. makes a decision and implements it) first, anticipating the follower's optimal response to its decision. Subsequently, the follower optimizes its own objective given the leader's action. The Stackelberg equilibrium problem can be formulated as a bilevel optimization problem. This hierarchical framework is applied in supply chain management,  providing an explicit optimization-based representation of this sequential decision-making process \citep{caselli2024bilevel,avraamidou2017multiparametric}. While two-player Stackelberg games are a primary application for bilevel optimization, bilevel optimization can be more broadly applied for the mathematical formulation of problems outside game theory, such as hyperparameter tuning in machine learning \citep{mackay2019self}.

Table~\ref{tab:carbon_tax_circularity} summarizes relevant studies involving multiple stakeholders modeled via Stackelberg games, with an emphasis on government policies such as carbon taxes and subsidies. These studies highlight the importance of using hierarchical approaches to model stakeholder interactions, capturing the sequential nature of decision-making and the interplay between environmental and economic objectives.

While several studies analyze taxes or incentives independently, and others consider fixed policy budgets, there is no study that simultaneously optimizes both instruments under a binding budget constraint within a bilevel policy-design framework. Considering joint policies helps capture the policy trade-offs and synergies that arise when a government must allocate limited financial resources across multiple incentive mechanisms. The resulting bilevel optimization framework provides novel theoretical and managerial insights, including (i) how subsidies can be coupled with carbon taxes to achieve circularity objectives without compromising industry profitability, (ii) how different policy mixes affect supply chain configurations and emissions outcomes in ways that cannot be inferred from tax-only or subsidy-only models, (iii) how budget limits restrict the feasible policy space and shape optimal technology adoption. The underlying intuition suggests that an increase in the carbon tax rate would reduce the corresponding incentive rate, while a decrease in the tax rate could enhance it. However, in real-world settings, this relationship is often complex and challenging to quantify, as it may exhibit non-linear dynamics arising from technological adoption decisions and policy interactions. To address this gap, the present work analyzes both policy mechanisms simultaneously within a circular supply chain context, considering the effect of different governmental budget levels to promote more circular and less-emitting supply chains.

\begin{table}
\centering
\footnotesize
\caption{Summary of studies on carbon tax and incentives impact on circular supply chains using Stackelberg game}
\label{tab:carbon_tax_circularity}

\resizebox{\textwidth}{!}{
\begin{tabular}{ >{\centering\arraybackslash}m{0.65cm} 
                 >{\centering\arraybackslash}m{3.6cm} 
                 >{\centering\arraybackslash}m{2.5cm} 
                 >{\centering\arraybackslash}m{3cm} 
                 >{\centering\arraybackslash}m{4.4cm} }
\hline
\textbf{Ref.} & \textbf{Players} & \textbf{Main Topic} & \textbf{Solution Approach} & \textbf{Contribution / Novelty} \\
\hline
\citep{hong2017optimizing} & Government, firms & Regional economics & Hybrid method (PDP + binary search + GA) & Government determines emission reduction targets via cap-and-trade \\
\citep{luo2022evaluating} & Manufacturer, retailer & Closed-loop supply chains & Analytical and numerical analysis & Carbon tax promotes remanufacturing but risks hindering circular activities if not carefully calibrated \\
\citep{liao2023optimal} & Manufacturer, recycler & Supply chain with recycling & Analytical and numerical analysis & Models government recycling incentives and tax effects on recycling supply chain \\
\citep{martelli2020optimization} & Government, multi-energy system owner & Multi-energy system & Nested metaheuristic (PGS-COM + CPLEX) & Considers subsidies and carbon tax simultaneously \\
\citep{Tsou2012} & Government, enterprise & Collection and recovery (recycling) stage & Karush–Kuhn– Tucker (KKT) approach & Develop subsidy and penalty strategies to minimize carbon emissions and costs in the recycling sector\\
\citep{chalmardi2019bi} & Government, supply chain owner & Sustainable supply chain network & Nested metaheuristic (SA + CPLEX) & Consider government incentives for circular supply chain design\\
\citep{camacho2023hierarchized} & Planning company, distributor & Sustainable supply chain network & Multi-start heuristic algorithm & Capturing economic-environmental Pareto fronts via bi-objective bilevel optimization.\\
\citep{ma2018incentive} & Supplier, manufacturer in notebook product line & Green product line & Nested genetic algorithm & Incentive-based strategy promoting sustainability \\
\citep{he2023hierarchical} & Government, multi-energy system operator & Multi-energy system & Teaching-learning based optimization & Lifecycle carbon reduction subsidy to promote sustainable energy adoption \\
\citep{mesrzade2023bilevel} & Government, supply chain operator & Product-specific supply chain & Grid search and CPLEX & Quantifies carbon tax impact on profit \\
\citep{rahmani2024competitive} & Government, supply chain operator & Supply chain for a pseudo product & Hybrid method (Quantum binary PSO + Benders decomposition) & Models government incentives as leader decisions in supply chain \\
This study & Government, industry & Packaging Waste Management & DOMINO framework (metaheuristic+analytical)&Jointly evaluates taxes and subsidies under budget constraints promoting circularity and emissions reduction\\
\hline
\end{tabular}}
\end{table}

\section{Methodology}

\subsection{Problem formulation - Bilevel optimization problem}

Bilevel programming problems involve a hierarchical structure where the lower-level optimization problem is embedded as a constraint within the upper-level problem. A general mixed-integer bilevel programming problem (B-MIP) is formulated as in Eq.~\ref{eqn:bilevel} \cite{sinha2017review}, featuring two objective functions: the upper-level objective $F$ and the lower-level objective $f$. In this work we formulate the policy making processes along with the response of the industry to the policy measures through bilevel optimization, where the upper-level represents the governmental decision making process, while the lower-level represents the response of product's manufacturers. Here, $\mathbf{x}$ denotes the vector of continuous variables, $\mathbf{y}$ denotes the vector of integer variables, and subscripts $u$ and $l$ refer to upper- and lower-level variables, respectively. Feasibility requires satisfaction of both upper- and lower-level constraints ($G$, $g$) and that the lower-level problem attains optimality for fixed upper-level variables.

\begin{singlespace}
\begin{align}
    \min_{\mathbf{x_u}, \mathbf{y_u}, \mathbf{x_l}, \mathbf{y_l}} & \quad 
    F(\mathbf{x_u}, \mathbf{y_u}, \mathbf{x_l}, \mathbf{y_l}) \notag \\
    \text{s.t.} & \quad 
    G(\mathbf{x_u}, \mathbf{y_u}, \mathbf{x_l}, \mathbf{y_l}) \leq 0, \notag \\
    & \quad 
    \mathbf{x_{l}}, \mathbf{y_{l}} \in \arg\min_{\mathbf{x_{l}}, \mathbf{y_{l}}} \big\{ f(\mathbf{x_u}, \mathbf{y_u}, \mathbf{x_l}, \mathbf{y_l}):\notag \\
    & \quad \hspace{4cm}
    g(\mathbf{x_u}, \mathbf{y_u}, \mathbf{x_l}, \mathbf{y_l}) \leq 0\big\}, \notag \\
    & \quad 
    \mathbf{x} = 
    \begin{bmatrix}
        \mathbf{x_{u}} \\
        \mathbf{x_{l}}
    \end{bmatrix}^T, \quad
    \mathbf{y} = 
    \begin{bmatrix}
        \mathbf{y_{u}} \\
        \mathbf{y_{l}}
    \end{bmatrix}^T, \notag \\
    & \quad 
    \mathbf{x} \in \mathbb{R}^n, \quad \mathbf{y} \in \mathbb{Z}^m.
    \label{eqn:bilevel}
\end{align}
\end{singlespace}

\begin{figure}
    \centering
    \includegraphics[width=0.5\linewidth]{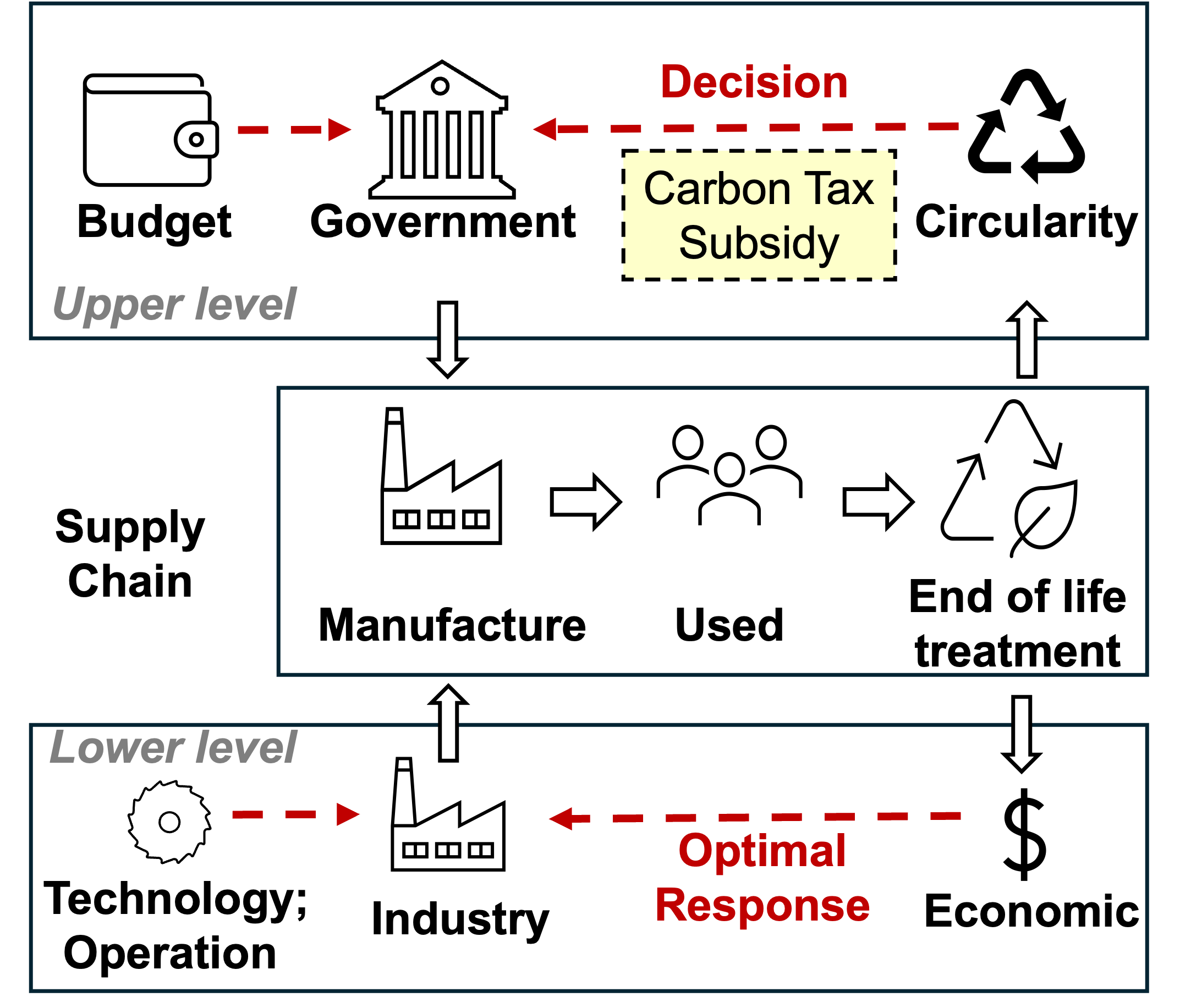}
    \caption{Bilevel circular supply chain model illustrating the hierarchical interaction between upper- and lower-level decision makers}
    \label{fig: Bilevel Circular supply chain}
\end{figure}

An example of hierarchical bilevel optimization for a circular supply chain is shown in Fig.~\ref{fig: Bilevel Circular supply chain} and Eq.~\ref{eqn: bilevel supplychain}. The supply chain involves two decision-makers: the (Regional) Government as the upper-level player and the Industry agents as the lower-level player. 

The government's objectives may encompass economic, environmental, and social sustainability considerations, guiding its policy instruments, which can include taxation, subsidies, and regulations. In this study, the government’s objective is explicitly defined in terms of environmental performance, but could be replaced or complemented by other economic or social criteria. The government aims to minimize the environmental impact by selecting the most sustainable technologies while balancing its budget constraints through taxes and subsidies.

At the lower level, the manufacturing industry seeks to maximize profit by optimizing technology choices in response to governmental policies, taxes, subsidies, and subject to operational constraints such as mass balance \& technology feasibility, and cost \& environmental impact calculations. The industry's actions can involve selecting manufacturing solutions, building infrastructure, and updating or decommissioning plants.

\begin{singlespace}
\begin{align}
\min_{taxes, subsidies} \quad & \text{Environmental Impact}(technologies) \notag\\
\text{s.t.} \quad & \text{Government Budget} \notag\\
\min_{technologies} \quad & \text{Cost}(taxes, subsidies, technologies) \notag\\
\text{s.t.} \quad & \text{Mass Balance and Technology Constraints} \notag\\
& \text{Cost and Environmental Impact Equation}
\label{eqn: bilevel supplychain}
\end{align}
\end{singlespace}

The interplay between these decision layers is modeled as a bilevel mixed-integer optimization problem where the government's environmental and economic goals form the upper-level problem, and the industry's profit maximization subject to technical and financial constraints constitutes the lower-level problem. This framework captures the conflict and interdependance between regulatory policies and industrial responses, facilitating integrated decision-making for sustainable development. The proposed bilevel supply chain model is then used to address the following research questions.

\begin{enumerate}
    \item What is the optimal technology to minimize GHG emissions considering both stakeholders?
    \item What are the corresponding carbon tax and subsidy levels for selecting the optimal technology?
    \item How does the budget influence carbon tax and subsidies?
    \item Given a carbon tax, how can the budget and subsidies be determined to select the optimal technology?
\end{enumerate}

\subsection{Mathematical model formulation and performance criteria}
In this study, two performance criteria are considered to quantify the environmental impact: the global warming potential (GWP), and the circularity index. The GWP quantifies the total GHG emissions of the supply chain design, while the circularity index is a composite index involving multiple indicators (total energy consumed, water withdrawal, water recycled or reused, hazardous and non-hazardous waste generated, etc) designed to quantify supply chain circularity. The circularity index used is obtained through the MICRON framework \citep{Baratsas2022}, and it ranges from 0 to 2, where values closer to 2 indicate a more circular system, and values closer to 0 reflect a more linear one. The economic impact is measured by minimizing the overall supply chain cost including production, waste management and transportation.

For the case of the government's objective being to minimize emissions, the resulting mathematical model incorporates both this environmental goal (as the upper-level objective) and the industry's economic objective (as the lower-level objective), alongside the effects of potential taxes and subsidies. The bilevel formulation is presented as follows:

\begin{singlespace}
\begin{align}
\min_{tax, sub} \quad &  GHG \  Emissions (p,u)  \notag\\
&=\sum_{p} \text{EP}_{p} + \sum_{p} \text{EPT}_{p} + \sum_{p,u} \text{ER}_{p,u} + \sum_{p,u} \text{ERT}_{p,u}
\quad \forall p \in P, u \in U \notag\\
\text{s.t.} \quad & Government \ Budget+Tax \ income \leq Sub \notag\\
\min_{p,u} \quad & Cost (p,u,tax,sub)\notag\\
&= \sum_p \text{CP}_{p} + \sum_{p} \text{CTP}_{p} -\sum_{p,u} \text{D}_{p, u} +tax \cdot  GHG \ Emissions - Sub \\
& \quad \quad \forall p \in P, u \in U \notag\notag\\
\text{s.t.} &\text{g}_M(tech,sub,tax) \leq 0  \ (Mass \ Balances \ and \ Technology \ Constraints)\notag\\
&\text{g}_P(tech,sub,tax) \leq 0  \ (Cost \ and\ Emissions\ Equations) \notag
\label{eqn: bilevel supplychain detail1}
\end{align}
\end{singlespace}

The upper level seeks to minimize the GHG emissions resulting from the selection of product $p$ and waste management technologies $u$ in the supply chain. These GHG emissions consider the impacts throughout the product life cycle, including contributions from production ($\text{EP}_p$) of product $p$, transportation ($\text{EPT}_p$), recycling ($\text{ER}_{p,u}$), and recycling transportation ($\text{ERT}_{p,u}$). To achieve this objective, the government controls two policy instruments: a carbon tax on GHG emissions and a subsidy for green technologies. These decision variables, 
carbon tax rate ($tax$) and subsidy rate ($sub$), are used to influence the decisions made by the industry in the lower level. The government's optimization problem is subject to a budget constraint, ensuring that the total amount allocated to subsidies does not exceed the available budget and tax revenue.

The lower-level problem aims to minimize the total cost which involves the product $p$ cost ($\text{CP}_p$), production transportation ($\text{CTP}_p$), the difference between the waste management cost and product sales ($D_{p,u}$) and the policies cost/benefit including the carbon tax cost given by the product of tax rate and the total GHG Emissions. Therefore, the cost depends on the product and technology selection, the carbon tax rate, and the subsidies decided by the government. The lower-level problem is subject to technological, budgetary, and process-related constraints. 

To encourage specific recycling practices, the government can adjust the specific subsidy rate ($sub_{p,t,r,a}$) of chosen routes. The total subsidy ($Sub$) is the sum of subsidies across all recycling routes, calculated by multiplying the specific subsidy rate by the number of bags for each route ($Bag_{p,t,r,a}$):

\begin{singlespace}
\begin{align}
Sub = \sum_{p,t,r,a}sub_{p,t,r,a}\cdot Bag_{p,t,r,a}
\end{align}
\end{singlespace}

When the government's objective is to maximize circularity, then the upper level changes in the following formulation: 

\begin{singlespace}
\begin{align}
\max_{tax, sub} \quad &  Circularity \ index (p,u)\notag\\
&=  \sum_{p} CEP_p + \sum_{p, \text{e}} CER_{p,\text{e}} 
\quad \forall p \in P, e \in E \notag\\
\text{s.t.} \quad & Government \ Budget+Tax \ income \leq Sub \notag\\
\min_{p,u} \quad & Cost (p,u,tax,sub)= \notag\\
&\sum_p \text{CP}_{p} + \sum_{p} \text{CTP}_{p} -\sum_{p,u} \text{D}_{p, u} +tax \cdot GHG \ Emissions - Sub \\
&  \quad \forall p \in P, u \in U \notag\\
\text{s.t.}  &\text{g}_M(tech,sub,tax) \leq 0  \ (Mass \ Balances \ and \ Technology \ Constraints)\notag\\
&\text{g}_P(tech,sub,tax) \leq 0  \ (Cost \ and\ Circularity\ Equations) \notag
\label{eqn: bilevel supplychain detail2}
\end{align}
\end{singlespace}

Here, the overall circularity index is given by the circularity calculated for the different product supply chain stages. In this case, it includes the circularity for the production stage (including transportation) $\text{CEP}_p$ and the waste management stage $\text{CER}_{p,e}$ for the product $p$ and the waste management scenarios $e$.

These bilevel formulations capture the interaction between government policy and industry behavior, where the government aims to reduce emissions or maximize circularity by imposing policies, and the industry reacts by optimizing its costs under these regulatory and technological limitations. A more detailed version of the model can be found in the Supplementary file.

\subsection{Solving bilevel optimization using the DOMINO framework}

Solving bilevel optimization problems is very challenging \citep{hansen1992new}. To solve the problem analytically, numerous approaches have been proposed. Karush–Kuhn–Tucker (KKT) reformulation is the most widely adopted, which can transform a linear bilevel optimization problem into a nonlinear optimization problem \cite{fortuny1981representation}, but cannot solve problems containing lower-level integers. B-MIP with both levels containing integers are even more challenging to solve \citep{kleinert2021survey}, prompting the development of different solution methods, such as Bender’s decomposition \citep{saharidis2009resolution}, multi-parametric optimization \citep{koppe2010parametric, avraamidou2019multi}, penalty functions \citep{dempe2005discrete}, and branch-and-bound/cut algorithms \citep{fischetti2018use}.

Analytical methods are effective for small-scale bilevel problems but are limited by scalability and computational demands for larger instances. For medium- to large-scale bilevel problems, hierarchical nested bilevel optimization frameworks using metaheuristic optimization are often used \citep{li2006hierarchical}. For a comprehensive review of this topic, see Camacho-Vallejo et al \citep{camacho2024metaheuristics}. 

Methods for solving bilevel optimization combining metaheuristic algorithms at the upper-level with exact mixed-integer linear programming (MILP) solvers at the lower level, such as the DOMINO framework \citep{beykal2020domino}, have been proposed. These approaches do not yield analytical solutions but produce bilevel-feasible approximations close to true optima by reformulating the bilevel problem as a single-level grey-box optimization. The procedure involves three main steps: (1) fixing the upper-level variables and solving the lower-level problems to optimality using exact solvers to evaluate the upper-level objective; (2) employing a metaheuristic solver to select new upper-level variable values based on the obtained objective; and (3) iterating step (1) and (2) until convergence or termination of the upper-level metaheuristic algorithm. 

Following the DOMINO framework, in this study, particle swarm optimization (PSO) \citep{kennedy1995particle} is employed as the upper-level metaheuristic method. The implementation is based on the work of Blank et al.\ \citep{pymoo}, using a swarm size of 10 and 200 iterations. It is important to note that metaheuristic methods do not guaranty a globally optimal solution, and both the results' quality and convergence speed are highly dependent on the choice of the starting point. In this study, we utilized domain knowledge to select different initialization points and report the best solutions found. While black-box global optimizers, such as the ISRES within the DOMINO framework, can approximate global optima, they entail significant computational overhead and have poor scalability as problem dimensionality increases \citep{beykal2020domino}. Consequently, such derivative-free approaches were excluded from this study. For obtaining an analytical solution to the MINLP-MILP bilevel optimization problem, one may consider reformulation-and-decomposition approaches combined with a nonlinear optimization solver \citep{yue2017stackelberg}. Nevertheless, there remains a critical need for developing algorithms that do not depend on intuition and domain expertise, that provide more generalizable solution frameworks.

\section{Coffee packaging supply chain case study}
 
Implementing circularity throughout the supply chain is essential for reducing environmental impacts and preserving resources. While much attention is often given to sustainable sourcing and production, waste management is a critical yet sometimes overlooked component of circular supply chains. Finding effective strategies among different waste management alternatives such as recycling and reuse, are necessary to close the materials loop and ensure that valuable resources are reintegrated into the production cycle rather than lost to landfills or incineration. Similarly, selecting products with lower environmental footprints and higher recyclability can enhance overall sustainability across the supply chain.

\begin{figure}
    \centering
    \includegraphics[width=0.5\linewidth]{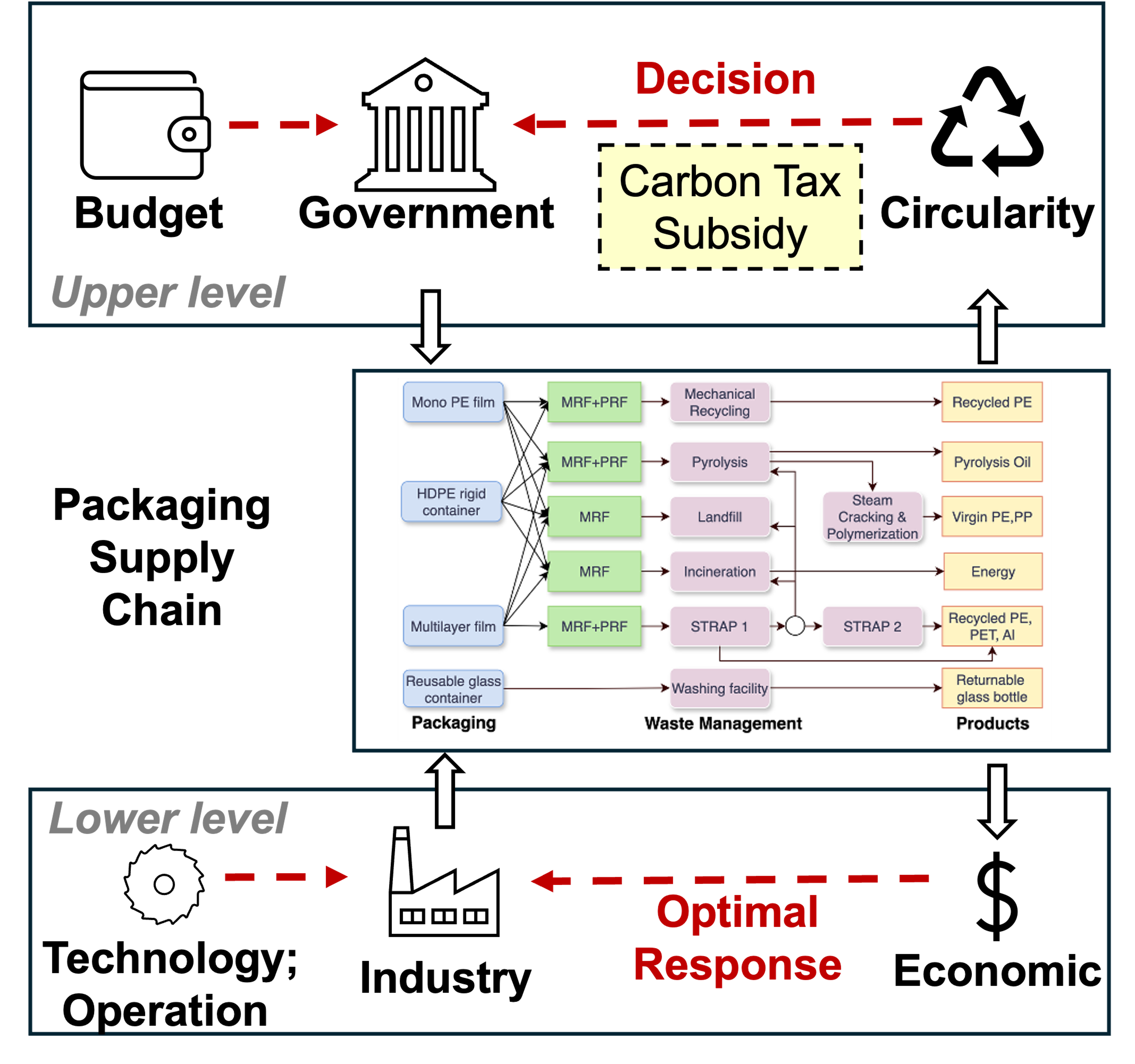}
    \caption{Case study representation of a bilevel circular packaging supply chain, highlighting the hierarchical interaction between the government and industry.}
    \label{fig: Bilevel Circular packaging supply chain}
\end{figure}

Collaboration among stakeholders is crucial for advancing circularity in the plastic packaging sector, as emphasized by Stumpf et al. \citep{stumpf2023circular}. Accordingly, incorporating policy instruments such as taxes and subsidies—and understanding their influence on stakeholders' behavior—is vital for driving sustainable outcomes. This case study builds upon the coffee packaging supply chain developed by Muñoz-Briones et al. \citep{munoz2025integrated}, and is adapted to incorporate bilevel optimization for modeling the hierarchical decision-making process between (regional) government and packaging manufacturers, balancing economic and environmental goals, as shown  in Fig.~\ref{fig: Bilevel Circular packaging supply chain}. 

The case study evaluated the economic and environmental performance of various packaging and waste management technologies to identify the optimal combination. The proposed framework integrates a MILP model with life cycle assessment (LCA), techno-economic analysis (TEA), and circularity assessment through the MICRON framework \citep{Baratsas2022}. An Extended Producer Responsibility (EPR) scenario is considered where the coffee production companies are responsible for the waste management of their product's packaging. The scope of the study includes the packaging production, waste management, and the corresponding transportation. Packaging options considered include multilayer and monolayer plastic films, rigid plastic containers, and reusable glass containers. The analysis covers a broad range of traditional and emerging waste management technologies, including landfill, incineration, glass washing, mechanical recycling, pyrolysis, and  solvent-based plastic recycling (STRAP™: Solvent-Targeted Recovery and Precipitation)\citep{Xu2025}. The case study considers a 1000 packaging units.

The upper and lower objectives reflect each stakeholder's interests, including cost minimization (manufacturer), GHG emission reduction, and promotion of circular economy ((regional) government). For the reformulation, the upper-level problem follows the environmental objective and controls the governmental policies, primarily focusing on tax and subsidy mechanisms designed to influence packaging material choices and waste management practices. The environmental performance is represented by two objectives: minimizing total GHG emissions and maximizing the overall circularity index. The lower-level problem follows the economic objective and models the supply chain decisions made by firms in response to these policies, including packaging and waste management technology selection. The economic objective is to minimize the increase in product price (cost), given by the difference between all the production, recycling and transportation costs and the revenue from product sales. 

 The proposed method is widely applicable to other industrial sectors. By changing the superstructure model's specific cost, emission, and circularity metrics, the bilevel optimization framework can capture the impact of taxes and incentives on industrial behavior. The results should be interpreted within the context of the specific case study analyzed, as they are contingent on the assumptions and data used in this analysis.

\section{Results and Discussion}
This study investigates two distinct upper-level governmental objectives: maximizing circularity and minimizing GHG emissions. We analyze the impacts of carbon tax, subsidies, and their combination on supply chain behavior under varying government budgets, industry and consumer behaviors (recycled glass loss rates), and transportation distances. Table~\ref{tab:Summary} summarizes the bilevel problems solved, their results in terms of packaging and recycling selected, along with notes on the policy efficacy towards achieving the objectives of the (regional) government. Table~\ref{tab:Summary} also includes the section numbers where more detailed discussion for each of the problem formulations and results can be found. 

\begin{table}[ht]
\centering
\caption{Technology selection of different waste management pathways under different policy handles.}

\begin{tabular}{%
    >{\centering\arraybackslash}m{2cm} 
    >{\centering\arraybackslash}m{2.4cm} 
    >{\centering\arraybackslash}m{5.3cm} 
    >{\centering\arraybackslash}m{2.8cm} 
    >{\centering\arraybackslash}m{1.8cm}} 
\hline
\textbf{Leader Objective} & \textbf{Policy Handles} & \textbf{Selected Pathway} & \textbf{Policy Efficacy} & \textbf{Section} \\
\hline
Min GHG Emission        & Carbon Tax         & multi-layer bag $\to$ STRAP™ $\to$ recycled PE/PET $\Longrightarrow$ multi-layer bag $\to$ landfill & Effective when Carbon Tax $>$\$4.3/kg-CO2 & ~\ref{sec: carbon taxes-GHG} \\
\hline
Max Circularity     & Carbon Tax         & multi-layer bag $\to$ STRAP™ $\to$ recycled PE/PET      & Not effective  & ~\ref{sec: carbon taxes-Circularity} \\ 
\hline
Min GHG Emission        & Subsidies          & multi-layer bag $\to$ STRAP™ $\to$ recycled PE/PET $\Longrightarrow$ multi-layer bag $\to$ landfill & Budget dependent & ~\ref{sec: incentives-GHG emissions} \\
\hline
Max Circularity     & Subsidies        & multi-layer bag $\to$ STRAP™ $\to$ recycled PE/PET $\Longrightarrow$ glass $\to$ glass washing & Budget dependent & ~\ref{sec: incentives-Circularity} \\
\hline
Min GHG Emission        & Subsidies and Carbon Tax  & multi-layer bag $\to$ STRAP™ $\to$ recycled PE/PET $\Longrightarrow$ glass $\to$ glass washing & Always effective & ~\ref{sec: Minimizing GHG tax-sub} \\
\hline
Max Circularity     & Subsidies and Carbon Tax  & multi-layer bag $\to$ STRAP™ $\to$ recycled PE/PET $\Longrightarrow$ multi-layer bag $\to$ landfill & Always effective & ~\ref{sec: max Circularity tax-sub}\\
\hline
\end{tabular}
\label{tab:Summary}
\end{table}

\subsection{Impact of Carbon Tax}
Carbon taxation is a widely employed instrument to promote GHG emissions reduction as it can incentivize businesses to lower their GHG emissions by making raw materials, processes, and products that are contributing to GHG emissions more costly and in effect promoting cleaner alternatives. Its popularity varies globally, as it is increasingly favored among policymakers and environmental advocates due to its economic efficiency, but policymakers and industry are concerned over its economic impacts and efficacy in reducing GHG emissions. Through the bilevel optimization formulation proposed in this work, we are able to explore the efficacy of carbon tax, and identify the level of carbon tax that will minimize the GHG emissions of the packaging industry case study. 

\subsubsection{Minimizing GHG Emissions}
\label{sec: carbon taxes-GHG}
The impact of the carbon tax in the supply chain where the governmental objective is minimizing GHG emissions, as shown in Eq.~\ref{eqn: carbon taxes-GHG}, exhibits a step change effect. As the tax rate increases, the supply chain abruptly shifts from a least cost scenario (multilayer bag that gets recycled using STRAP™ to produce recycled PE, Al, PET) to a least emission scenario (multilayer bag that goes to landfill). The sudden change occurs at a carbon tax rate threshold of \$4.3/kg-CO$_{2,\mathrm{eq}}$. 

This result indicates that any carbon tax rate bellow \$4.3/kg-CO$_{2,\mathrm{eq}}$ will result to an increase cost for the packaging supply chain, without leading to any technology change or GHG emission reduction. Notably, a high tax level is required given the lowest-cost recycling option (STRAP) is profitable and does not generate significantly higher emissions than landfilling.

\begin{singlespace}
\begin{align}
\min_{carbon \ tax} \quad & \text{GHG Emissions}(technologies) \notag\\
\text{s.t.} \quad & \text{Government Budget} \notag\\
\min_{technologies} \quad & \text{Cost}(carbon\ tax, technologies) \notag\\
\text{s.t.} \quad & \text{Mass \ Balance \  and \ Technology \ Constraints} \notag\\
& \text{Cost \ and \ Emissions \ Equations}
\end{align}
\label{eqn: carbon taxes-GHG}
\end{singlespace}

\subsubsection{Maximizing Circularity}
\label{sec: carbon taxes-Circularity}
While carbon taxation is a widely employed instrument to promote GHG emissions reduction and circularity in supply chains, its efficacy to driving a transition from linear economy to circular economy is not clear, with researchers suggesting that it cannot directly drive towards circularity, or suggesting that it needs complementary measures \citep{panza2025role}. The case study under consideration here shows that carbon taxation effectiveness for transition towards circularity is limited.

The carbon tax has no impact on the bilevel supply chain, when the government objective is to maximize circularity as defined in Eq.~\ref{eqn: carbon taxes-Circularity}. Therefore, increasing the tax rate does not alter the supply chain configuration, which remains in the least-cost scenario (Multilayer bag that goes to STRAP™ and produces recycled PE, Al, PET) rather than shifting to the most circular, lowest-emission configuration (Glass container that goes to glass washing). This is because, despite higher circularity, the glass washing scenario generates more GHG emissions than the least-cost alternative. Consequently, the overall cost difference between scenarios increases with the carbon tax, preventing the transition to the most circular solution through carbon taxation alone.

\begin{singlespace}
\begin{align}
\max_{carbon \ tax} \quad & \text{Circularity}(technologies) \notag\\
\text{s.t.} \quad & \text{Government Budget} \notag\\
\min_{technologies} \quad & \text{Cost}(carbon \ tax, technologies) \notag\\
\text{s.t.} \quad & \text{Mass \ Balance \  and \ Technology \ Constraints} \notag\\
& \text{Cost \ and \ Circularity \ Equations}
\label{eqn: carbon taxes-Circularity}
\end{align}
\end{singlespace}

\subsection{Impact of Subsidies}

In this section, subsidies are allocated to specific supply chain configurations to incentivize their selection. Given sufficient budget, the optimal governmental objectives of minimizing GHG emissions and maximizing circularity can be achieved. When the budget is limited, the objectives improve progressively as the subsidy budget increases. This differs from a carbon tax rate exhibiting a discontinuous transition at a specific threshold.

\subsubsection{Minimizing GHG Emissions}
\label{sec: incentives-GHG emissions}
\begin{figure}
    \centering
    \includegraphics[width=0.8\linewidth]{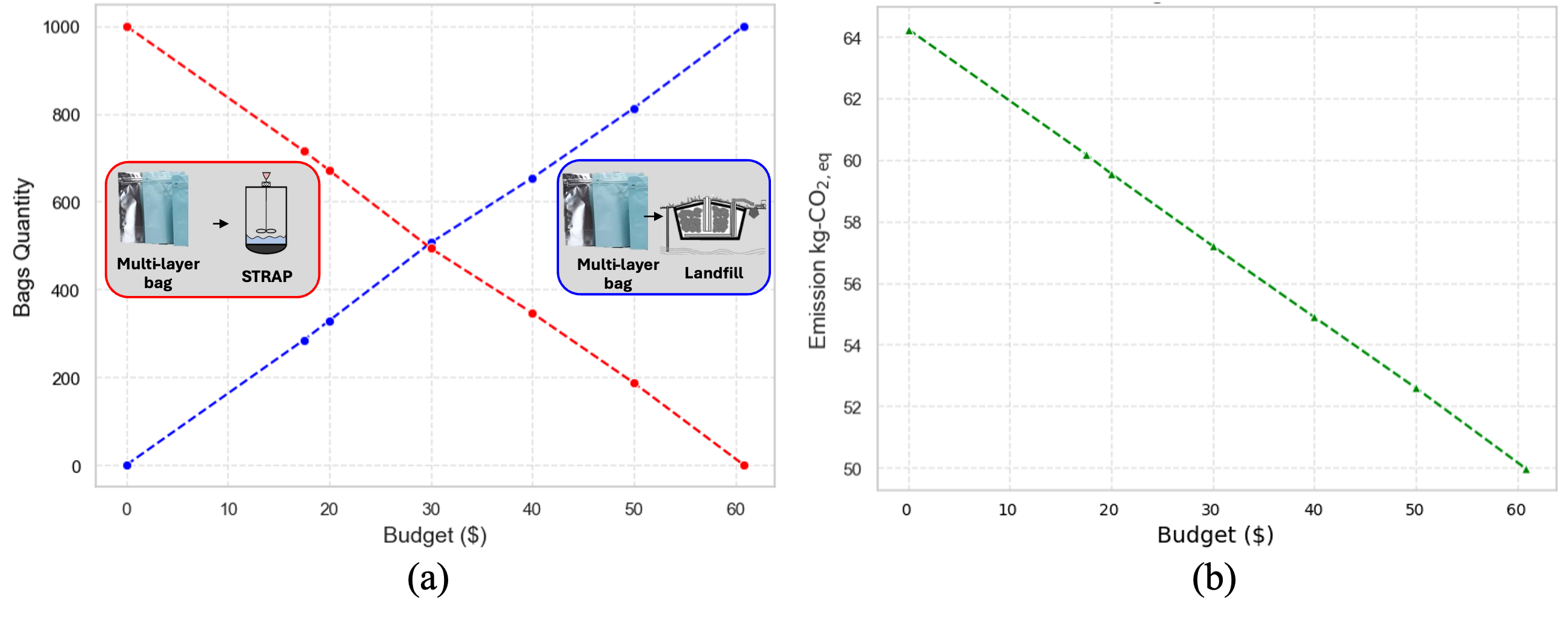}
    \caption{Impact of subsidies on minimizing GHG emissions: Relationship between government budget and (a) types and quantities of packaging end-of-life treatments, and (b) resulting GHG emission levels.}
    \label{fig:subsidy emission}
\end{figure}

The impacts of the incentives on the bilevel supply chain model where the government is minimizing the GHG emission is defined as Eq.~\ref{eqn: incentives-GHG emissions}, and the results are shown in Fig.~\ref{fig:subsidy emission}

\begin{singlespace}
\begin{align}
& \min_{Subsidies} \quad \textit{GHG Emissions}(\text{technologies}) \notag\\
& \text{s.t.} \quad \text{Government Budget} \geq \text{Subsidies} \notag\\
& \quad \min_{\text{technologies}} \quad \text{Cost}(\text{subsidies}, \text{technologies}) \notag\\
\text{s.t.} \quad & \text{Mass \ Balance \  and \ Technology \ Constraints} \notag\\
& \text{Cost \ and \ Emissions \ Equations}
\label{eqn: incentives-GHG emissions}
\end{align}
\end{singlespace}

The subsidies incentivize the configuration with the lowest GHG emissions, specifically the multilayer bag to landfill process (blue line in Fig.~\ref{fig:subsidy emission}(a)). When the budget is zero, there is no incentive to give out. Hence, the supply chain selects the least-cost solution (multilayer bag to STRAP™ to produce recycled PE, Al, PET). As the budget increases, the packaging manufacturer linearly shifts from the least-cost process toward the lowest-emission pathway. A subsidy rate of \$0.061 per multilayer bag enables the complete transition from STRAP™ to landfill. This change results in a total decline in GHG emissions as more bags are diverted to landfill, decreasing from 64.24 kg-CO$_{2,\mathrm{eq}}$ to 49.97 kg-CO$_{2,\mathrm{eq}}$, as shown in Fig.~\ref{fig:subsidy emission}(b). The supply chain selects the least emission solution in the end, where the budget exceeds \$61 for this 1000 bags scenario.

\subsubsection{Maximizing Circularity}
\label{sec: incentives-Circularity}
\begin{figure}
    \centering
    \includegraphics[width=0.8\linewidth]{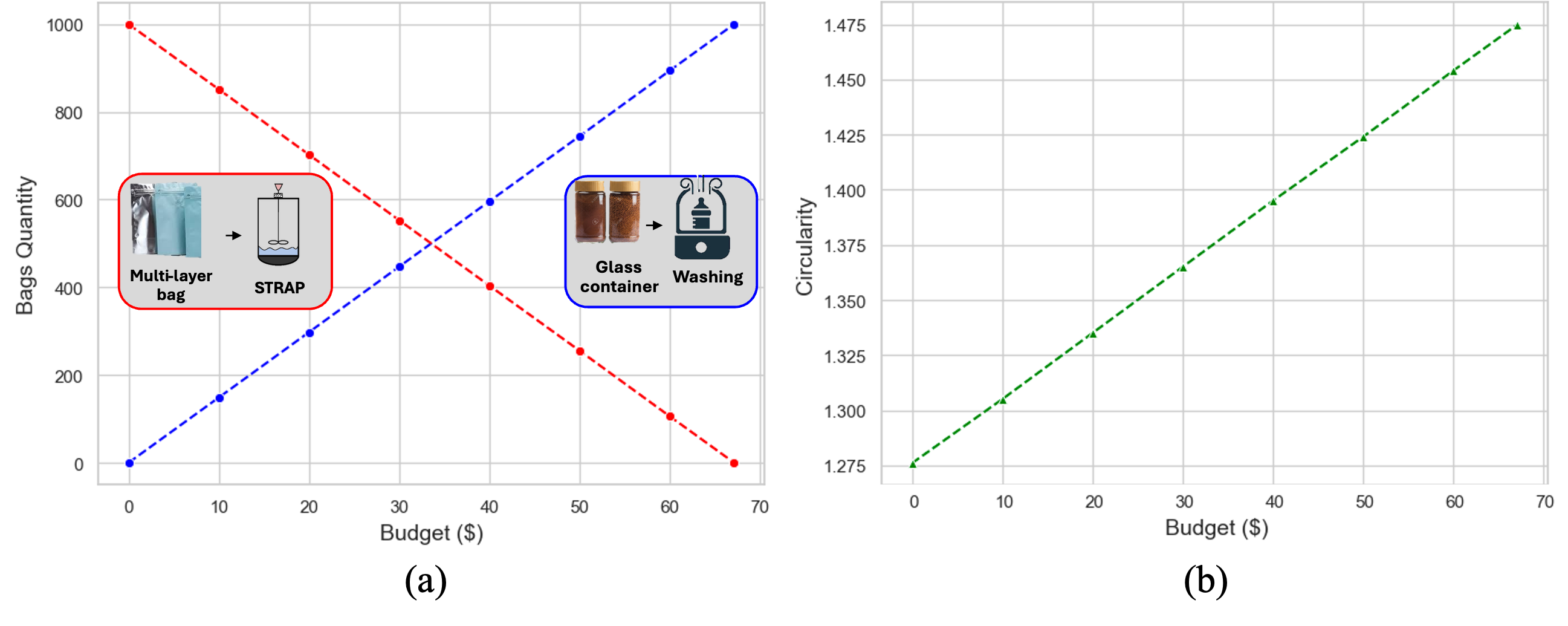}
    \caption{Impact of subsidies on maximizing circularity: Relationship between government budget and (a) packaging end-of-life treatments types and quantities, and (b) Circularity index levels.}
    \label{fig:subsidy circularity}
\end{figure}

The impacts of the subsidies on the bilevel supply chain model where the government is maximizing the overall circularity is defined as Eq.~\ref{eqn: incentives-Circularity}. The government subsidizes the selection of the most circular configuration (Glass container $\rightarrow$ Glass washing) at a rate of \$0.067 per glass jar. Initially, the supply chain selects the least-cost solution (multilayer bag to STRAP™ to produce recycled PE, Al, PET). When the budget increases, the packaging manufacturer shifts linearly from the least-cost process toward the highest circularity configuration, as shown in Fig.~\ref{fig:subsidy circularity}(a). As a result, the overall circularity increases from 1.275 to 1.475 (Fig.~\ref{fig:subsidy circularity}(b)) as a greater fraction of glass jars are directed to glass washing. Ultimately, when the budget exceeds \$67 for the 1000-bag scenario, the supply chain adopts the most circular solution.

\begin{singlespace}
\begin{align}
& \max_{Subsidies} \quad \textit{Circularity}(\text{technologies}) \notag\\
& \text{s.t.} \quad \text{Government Budget} \geq \text{Subsidies} \notag\\
& \quad \min_{\text{technologies}} \quad \text{Cost}(\text{subsidies}, \text{technologies}) \notag\\
\text{s.t.} \quad & \text{Mass \ Balance \  and \ Technology \ Constraints} \notag\\
& \text{Cost \ and \ Circularity \ Equations}
\label{eqn: incentives-Circularity}
\end{align}
\end{singlespace}

\subsection{Combination of Policy Handles: Impact of Carbon Tax and Subsidies}

This study examines the combined interactions of two government policy tools, carbon tax and subsidies within the supply chain, which remain underexplored in the literature. The analysis is based on a government budget constraint ($\text{Government Budget} + \text{Tax Income} \geq \text{Subsidies}$) \citep{milyani2018optimally}. This constraint relates budget and carbon tax revenue, to subsidy allocation.

\subsubsection{Minimizing GHG Emissions}
\label{sec: Minimizing GHG tax-sub}
First, a zero-budget scenario was analyzed to address the first two research questions, where the government aims to minimize emissions while applying both policy instruments Eq.~\ref{eqn: minimizing GHG tax-sub}. The results, presented in Table~\ref{tab:solution_budget_0}, show that the technology with the lowest emissions (Multilayer bag $\rightarrow$ Landfill) is selected. A carbon tax rate of \$0.9/kg-CO$_2$ generates sufficient revenue (\$47) to meet the zero budget constraint, providing incentives (\$47) for the selection of the least-emitting configuration while economically penalizing the other alternatives.

\begin{singlespace}
\begin{align}
\min_{\text{Carbon tax \& Subsidies}} \quad & \text{GHG Emissions}(\text{technologies}) \notag\\
\text{s.t.} \quad & \text{Government Budget} + \text{Tax Income} \geq \text{Subsidies}\notag\\
\min_{\text{technologies}} \quad & \text{Cost}(\text{carbon tax}, \text{subsidies}, \text{technologies}) \notag\\
\text{s.t.} \quad & \text{Mass \ Balance \  and \ Technology \ Constraints} \notag\\
& \text{Cost \ and \ Emissions \ Equations}
\label{eqn: minimizing GHG tax-sub}
\end{align}
\end{singlespace}


\begin{table}[ht]
\centering
\caption{Optimization results for the three different objectives for a zero-budget scenario.}
\begin{tabularx}{\textwidth}{l>{\centering\arraybackslash}p{3cm}>{\centering\arraybackslash}p{3cm}>{\centering\arraybackslash}X}
\hline
& \makecell{\textbf{Min GHG} \\ \textbf{Emission}} 
& \makecell{\textbf{Max} \\ \textbf{Circularity}} 
& \makecell{\textbf{Most} \\ \textbf{Profitable}} \\ \hline

\textbf{Technology selected} 
& Multilayer $\to$ Landfill 
& Glass $\to$ Glass washing 
& Multi-layer bag $\to$ STRAP™ $\to$ recycled PE/PET \\ 

\textbf{Carbon Emission} 
& 49.97 kg-CO$_2$ 
& 50.08 kg-CO$_2$ 
& 64.24 kg-CO$_2$ \\ 

\textbf{Cost}
& -\$61.52 
& -\$67.94 
& -\$0.93 \\ 

\textbf{Circularity Index }
& 1.18 
& 1.48 
& 1.28 \\ 

\textbf{Carbon tax rate} 
& \$0.9/kg-CO$_2$ 
& \$1.05/kg-CO$_2$ 
& --- \\ 

\textbf{Carbon tax} 
& \$47 
& \$52 
& --- \\ 

\textbf{Incentives rate} 
& \$0.047/container 
& \$0.052/container 
& --- \\ 

\textbf{Incentives} 
& \$47 
& \$52 
& --- \\ \hline
\end{tabularx}
\label{tab:solution_budget_0}
\end{table}

\begin{figure}
    \centering
    \includegraphics[width=0.8\linewidth]{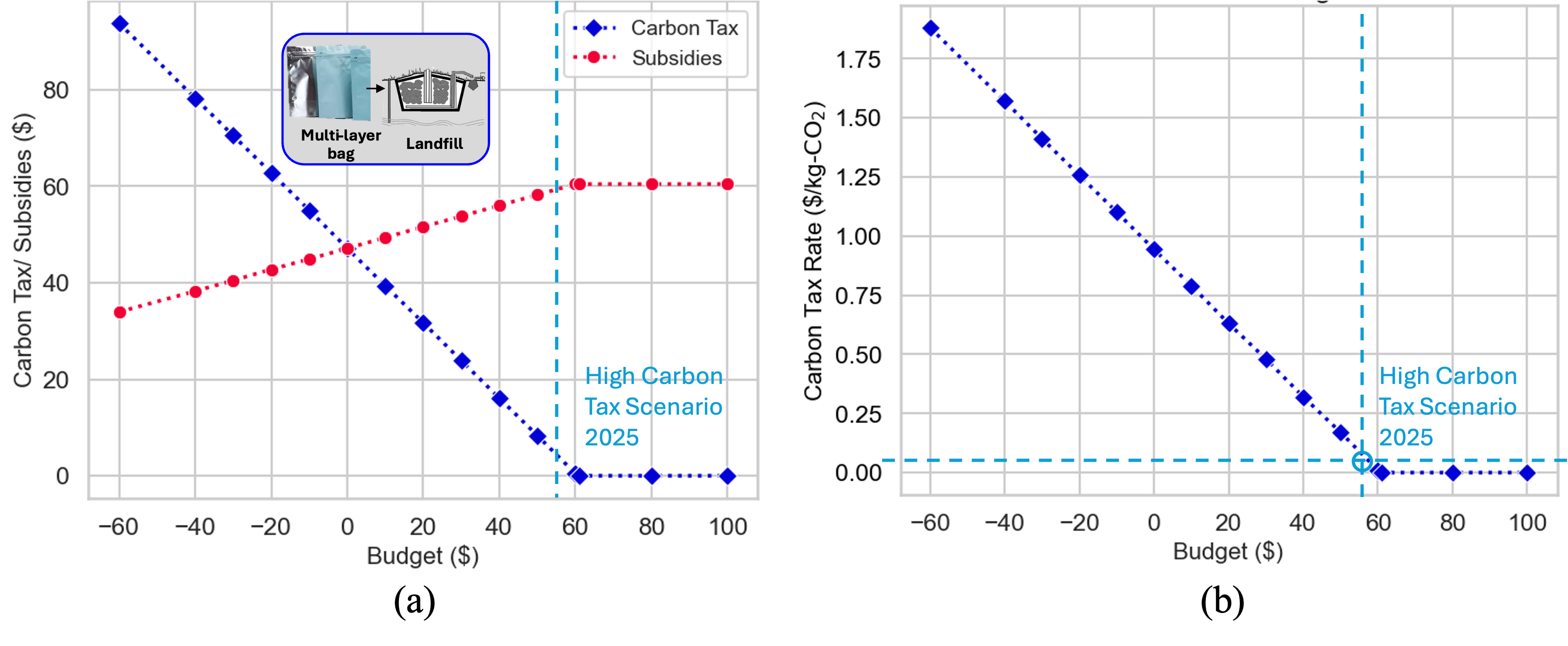}
    \caption{Impact of carbon tax and subsidies on minimizing GHG emissions: Relationship between government budget and (a) Carbon Tax revenue and subsidies, and (b) Carbon tax rate, highlighting a high carbon tax benchmark in 2025 at 0.1 \$/kg-CO$_{2,\mathrm{eq}}$.}
    \label{fig:subsidy tax emission}
\end{figure}

A sensitivity analysis of government budget ranges from \$-60 to \$100 is carried out to address question 3. Positive budget values indicate net government spending on the supply chain, and negative values represent net government revenue generation. Fig.~\ref{fig:subsidy tax emission}a shows the relationship between government budget, tax, and subsidies. This piecewise linear behavior clarifies how budget adjustments impact tax and subsidy dynamics in the supply chain. For budgets below \$61, the tax (blue dotted line) decreases linearly with increasing budget, while subsidies (red dotted line) increase linearly.  For budgets above \$61, the carbon tax remains zero, and subsidies remain unchanged. All government spending flows directly to the coffee company, and the net profit increase corresponds to the budget, subsidy minus carbon tax curve, resulting in a 45-degree line passing through the origin. 

To address question 4, a fixed carbon tax rate can be represented by a horizontal line, exemplified by the cyan dotted line at \$0.1/kg-CO$_2$ in Fig.~\ref{fig:subsidy tax emission}b. This rate corresponds to a high carbon tax rate benchmark globally in 2025. The intersection of this fixed tax rate and the carbon tax curve (blue dotted line) determines the required budget of approximately \$55, which is necessary to drive the supply chain to its minimum emission configuration. The associated carbon tax and subsidies, indicated by the intersection of the cyan dotted and tax/subsidies lines, are approximately \$4 in carbon tax and \$59 in subsidies, respectively.

\subsubsection{Maximizing Circularity}
\label{sec: max Circularity tax-sub}
\begin{figure}
    \centering
    \includegraphics[width=0.8\linewidth]{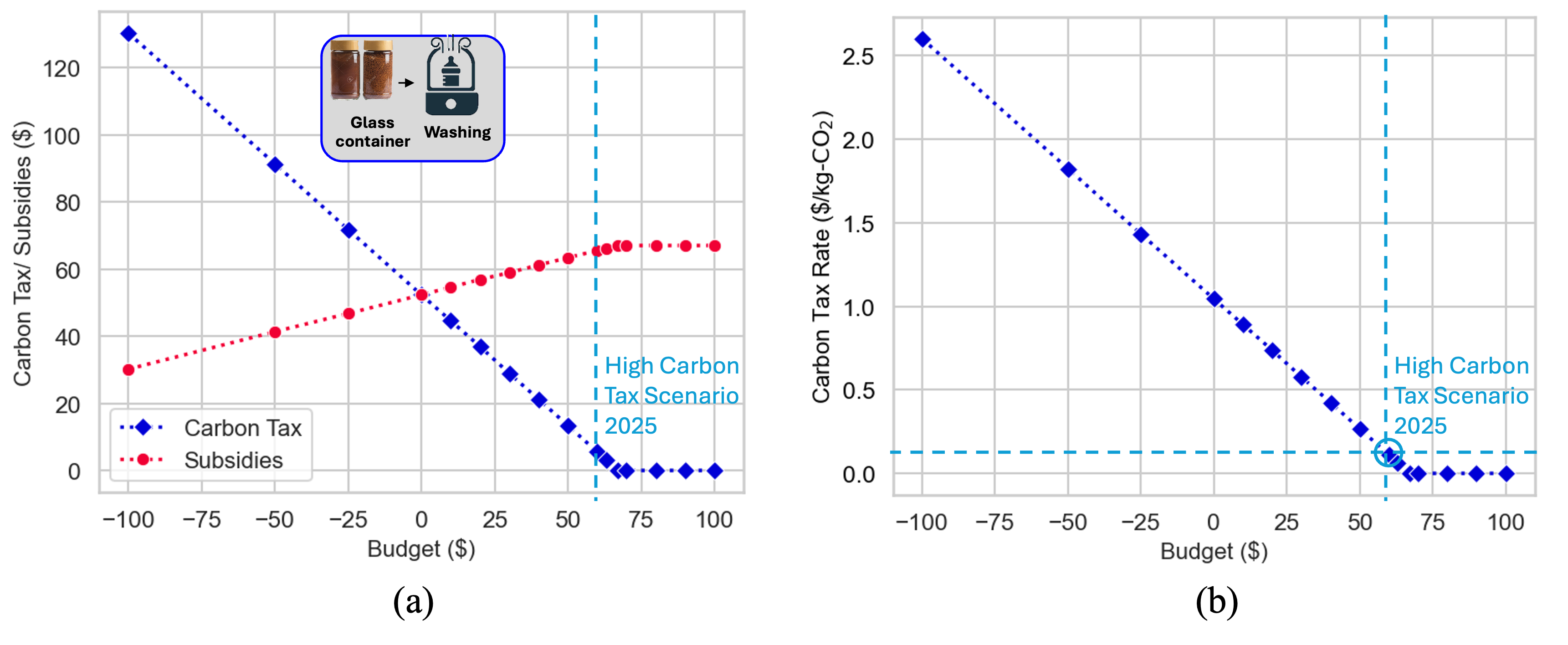}
    \caption{Impact of carbon tax and subsidies on maximizing circularity: Relationship between government budget and (a) Carbon Tax revenue and subsidies, and (b) Carbon tax rate, highlighting a high carbon tax benchmark in 2025 at 0.1 \$/kg-CO$_{2,\mathrm{eq}}$.}
    \label{fig:subsidy tax circularity}
\end{figure}

Fig.~\ref{fig:subsidy tax circularity} presents the results for maximizing circularity by solving Eq.~\ref{eqn:  max Circularity tax-sub}. The four questions outlined in Section~\ref{sec: Minimizing GHG tax-sub} are addressed similarly here. The results displayed in Table ~\ref{tab:solution_budget_0} show that the most circular scenario is selected (Glass container $\rightarrow$ Glass washing) when a carbon tax rate of \$1.05 per kg-CO$_2$eq is imposed to all the technologies and an incentive of \$0.052 per packaging unit is provided. The figure also illustrates the sensitivity of tax and subsidies relative to the budget. A scenario with a fixed carbon tax rate can also be evaluated. For instance, at a high carbon tax rate benchmark of \$0.1/kg-CO$_{2,\mathrm{eq}}$ (highlighted in the cyan dashed line), a budget of \$60 is required to provide the corresponding subsidies. 

\begin{singlespace}
\begin{align}
\max_{\text{Carbon tax \& Subsidies}} \quad & \text{Circularity}(\text{technologies}) \notag\\
\text{s.t.} \quad & \text{Government Budget} + \text{Tax Income} \geq \text{Subsidies}\notag\\
\min_{\text{technologies}} \quad & \text{Cost}(\text{carbon tax}, \text{subsidies}, \text{technologies}) \notag\\
\text{s.t.} \quad & \text{Mass \ Balance \  and \ Technology \ Constraints} \notag\\
& \text{Cost \ and \ Circularity \ Equations}
\label{eqn: max Circularity tax-sub}
\end{align}
\end{singlespace}

Sections~\ref{sec: Minimizing GHG tax-sub} and~\ref{sec: max Circularity tax-sub} quantify the trade-off between carbon tax rates and subsidy allocations under a fixed governmental budget. Despite imposing a relatively high carbon tax rate of \$0.1 per kg CO$_{2,\mathrm{eq}}$ in 2025, it is observed that the generated revenue covers only a minor fraction of the subsidies required for both GHG minimization and circularity maximization scenarios. Consequently, governmental budget allocation remains the primary driver for achieving these objectives. Among the analyzed configurations, the most economically favorable solution: recycling multi-layer bags via the STRAP™ method has low carbon intensity, characterized by small emissions combined with significant cost savings relative to alternative pathways. Therefore, shifting away from this economically favorable configuration through carbon taxation alone is challenging. As shown in Fig.~\ref{fig:subsidy tax emission} and~\ref{fig:subsidy tax circularity}, the solution under a high carbon tax closely coincides with that under no tax.

For future work, carbon tax cap is another aspect to consider that can lead to interesting behaviors in some cases. Current assumptions in sections~\ref{sec: Minimizing GHG tax-sub} and~\ref{sec: max Circularity tax-sub} have no upper bounds on the carbon tax limit. This assumption results in the consistent selection of the least-emission and maximum-circularity scenarios, since the government can provide unlimited subsidies by raising the carbon tax to infinite levels, even when the budget is zero or negative.

\subsection{Sensitivity Analysis}

In this section, we explore the sensitivity of two key parameters: the distance to the glass washing facility and the percentage of glass loss.

\subsubsection{Glass washing distance}
\begin{figure}
    \centering
    \includegraphics[width=0.8\linewidth]{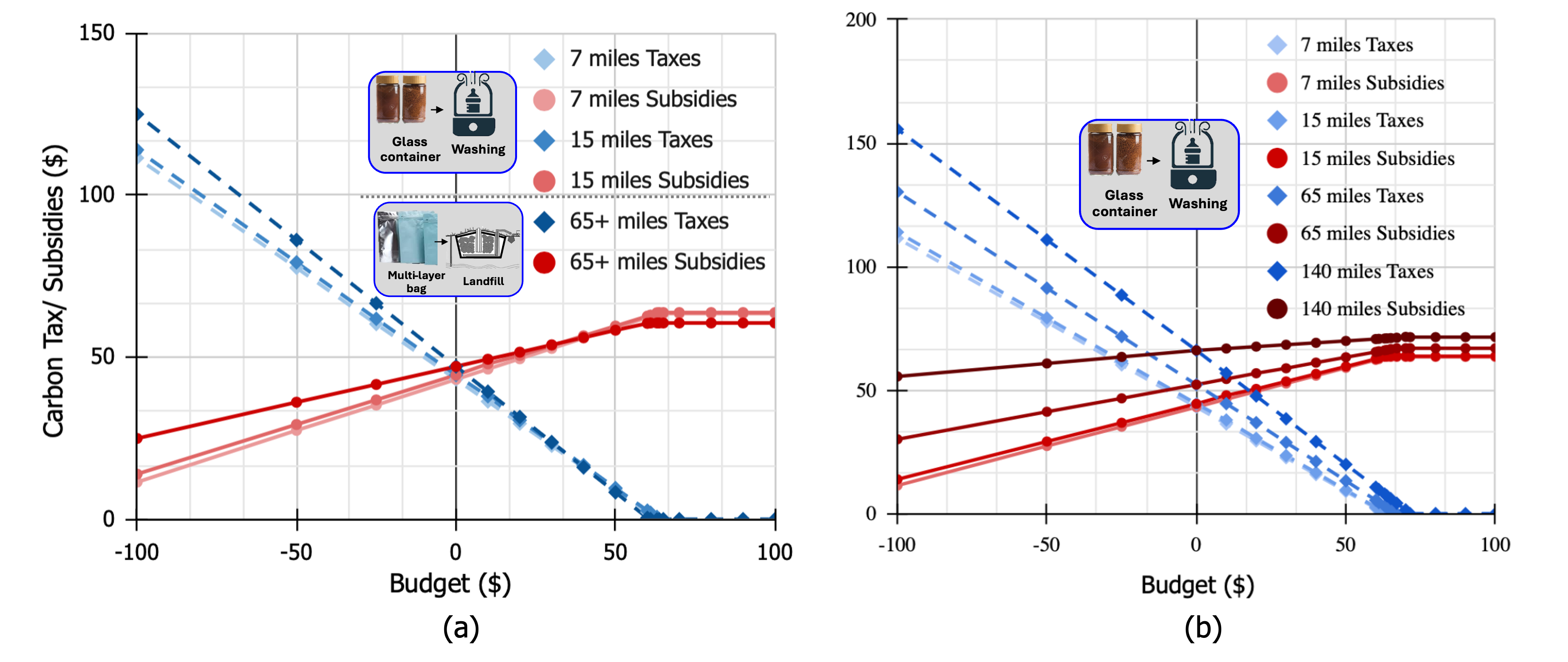}
    \caption{Sensitivity analysis of glass washing facility distance illustrating the relationship between carbon tax and subsidies as a function of budget: (a) minimizing GHG emissions for distances of 7, 15, and 65+ miles; (b) maximizing circularity for distances of 7, 15, 65, and 140 miles.}
    \label{fig:max Emi_Cir dis.png}
\end{figure}

The analysis specifically focuses on the glass washing distance, transportation from the collection center to glass washing facility, following the findings of Muñoz-Briones et al. \citep{munoz2025integrated}, which demonstrate that the most circular configuration does not necessarily align with the least-emitting alternative. Given that the glass packaging is much heavier than alternative packaging materials, the increased weight results in higher transportation costs and emissions, potentially influencing the selection of end-of-life treatment technologies. The analysis considers four different distances: 7, 15, 65, and 140 miles to represent different local and non-local distances. These distances are considered in this case study as well, as they can represent different alternatives like
a local collection facility, municipal facility, multi-city regional plant, or state-based centralized facility. The bilevel model considering both environmental metrics: minimizing GHG emissions and maximizing circularity, are evaluated and illustrated in Figs.~\ref{fig:max Emi_Cir dis.png}

When minimizing GHG emissions, the lowest emissions are achieved by employing the \emph{glass} $\rightarrow$ \emph{glass washing} technology for the 7- and 15-mile scenarios, while the \emph{multilayer bag} $\rightarrow$ \emph{landfill} option is favored for distances of 65 miles and above. The corresponding carbon tax and budget requirements are illustrated in Fig.~\ref{fig:max Emi_Cir dis.png}(a). 

The incentives and carbon tax necessary to offset emissions increase as the glass washing transport distance extends, which is expected as more emissions will be generated. The trend observed for distances of 65 miles and more differs from that of the lower distance cases (7 and 15 miles) and becomes less sensitive to changes in the budget. This occurs due to the adoption of a different technology. Since the 65-mile and above scenarios select a technology unrelated to the glass washing distance, the computed carbon tax and incentives lines remain the same. 

In the analysis of maximizing circularity, the highest circularity is achieved by employing the \emph{glass}~$\rightarrow$~\emph{glass washing} technology across all scenarios. The corresponding carbon tax and budget requirements are illustrated in Fig.~\ref{fig:max Emi_Cir dis.png}(b).

Overall, as the glass washing distance increases, and assuming no technological shift occurs, emissions increase, also resulting in a slight rise in the carbon tax rate to impose a higher penalty. Consequently, the overall carbon tax, defined as the product of these two factors, exhibits a steeper slope. Oppositely, due to budget constraints, the slope of the subsidies constraint decreases as the glass washing distance increases, as illustrated in both Fig~\ref{fig:max Emi_Cir dis.png}(a,b).

In Fig.~\ref{fig:max Emi_Cir dis.png}(a), where a technology shift occurs, the multilayer bag $\rightarrow$ landfill technology is less costly than the glass $\rightarrow$ glass washing technology; therefore, the required incentives are lower when there is a sufficient budget (greater than \$60). Moreover, we know that the \emph{multilayer bag} $\rightarrow$ \emph{landfill} pathway has a larger carbon emission than the 7- and 15-mile \emph{glass} $\rightarrow$ \emph{glass washing} pathways (otherwise, it would be the selected solution, which is not the case). Consequently, the slope of the carbon tax is larger, as more carbon tax is needed. However, when comparing the solutions selected for the 65-miles scenario for both environmental objectives, the carbon tax slope for multilayer film to landfill (Fig~\ref{fig:max Emi_Cir dis.png}(a): 65+ carbon tax) is less steep than the 65-mile \emph{glass} $\rightarrow$ \emph{glass washing} pathway (Fig.~\ref{fig:max Emi_Cir dis.png}(b):65-mile carbon tax), because the carbon emission is less when minimizing GHG emissions.

\subsection{Glass loss percent}

Another key parameter affecting sensitivity is the glass loss percentage, which represents the proportion of glass lost during the return and reuse cycle. This typically results from breakage, improper collection, or rejection when cleaning standards are not met. For the base case, we have chosen a loss rate of 3.13\%, representing the average operation of four bottling companies in northern Italy \citep{tua2020reusing}. In this section, we explore the impact of glass loss rates of 1\%, 3.13\%, and 10\% on the objectives of minimizing GHG emissions and maximizing circularity while fixing our glass washing facility distance as 65 miles. This analysis provides a comprehensive understanding of how the glass loss percentage influences the implementation of carbon tax and subsidy policies within the bilevel optimization framework.

\begin{figure}
    \centering
    \includegraphics[width=0.75\linewidth]{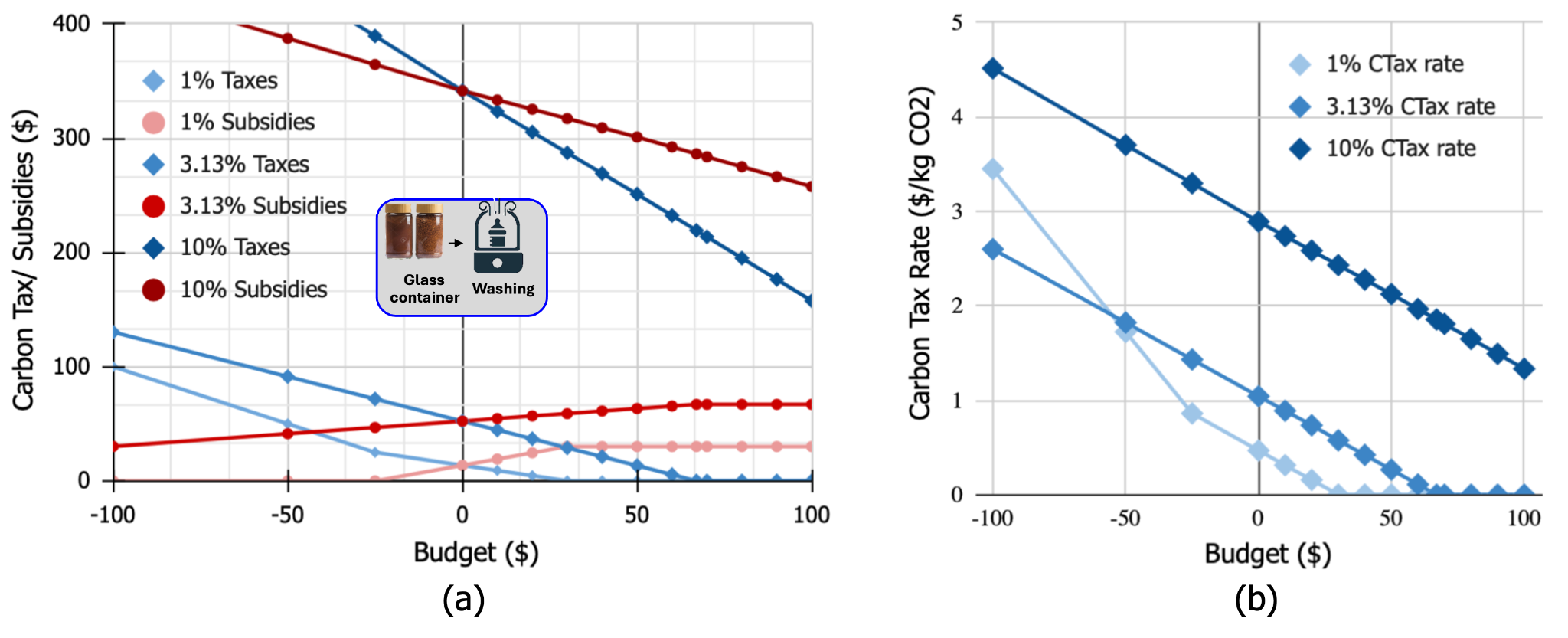}
    \caption{Sensitivity analysis of minimizing GHG emissions for glass loss rates of 1\%, and 3.13+\%. (a) Carbon tax revenue and subsidies as functions of government budget. (b) Carbon tax rate versus government budget.}
    \label{fig:min Emission gloss}
\end{figure}

When minimizing GHG emissions, the 1\% glass loss scenario selects the \emph{glass} $\rightarrow$ \emph{glass washing} pathway, while the 3.13\% glass loss scenario selects the \emph{multilayer bag} $\rightarrow$ \emph{landfill} pathway. As the glass loss percentage increases, GHG emissions also increase. Consequently, for all glass loss percentages above 3.13\%, the selected pathway remains constant (\emph{multilayer bag} $\rightarrow$ \emph{landfill} and no longer depends on the glass loss percentage. Notably, the 1\% glass loss also results in lower costs; therefore, the subsidies required are less than those for other pathways. Additionally, we observed a change in the slope of the carbon tax and carbon tax rate for 1\% glass loss scenario. This occurs because the solution enters a region (Budget$<$-25\$) where no subsidies are required to drive the technology change. The results are illustrated in Fig.~\ref{fig:min Emission gloss}.

\begin{figure}
    \centering
    \includegraphics[width=0.75\linewidth]{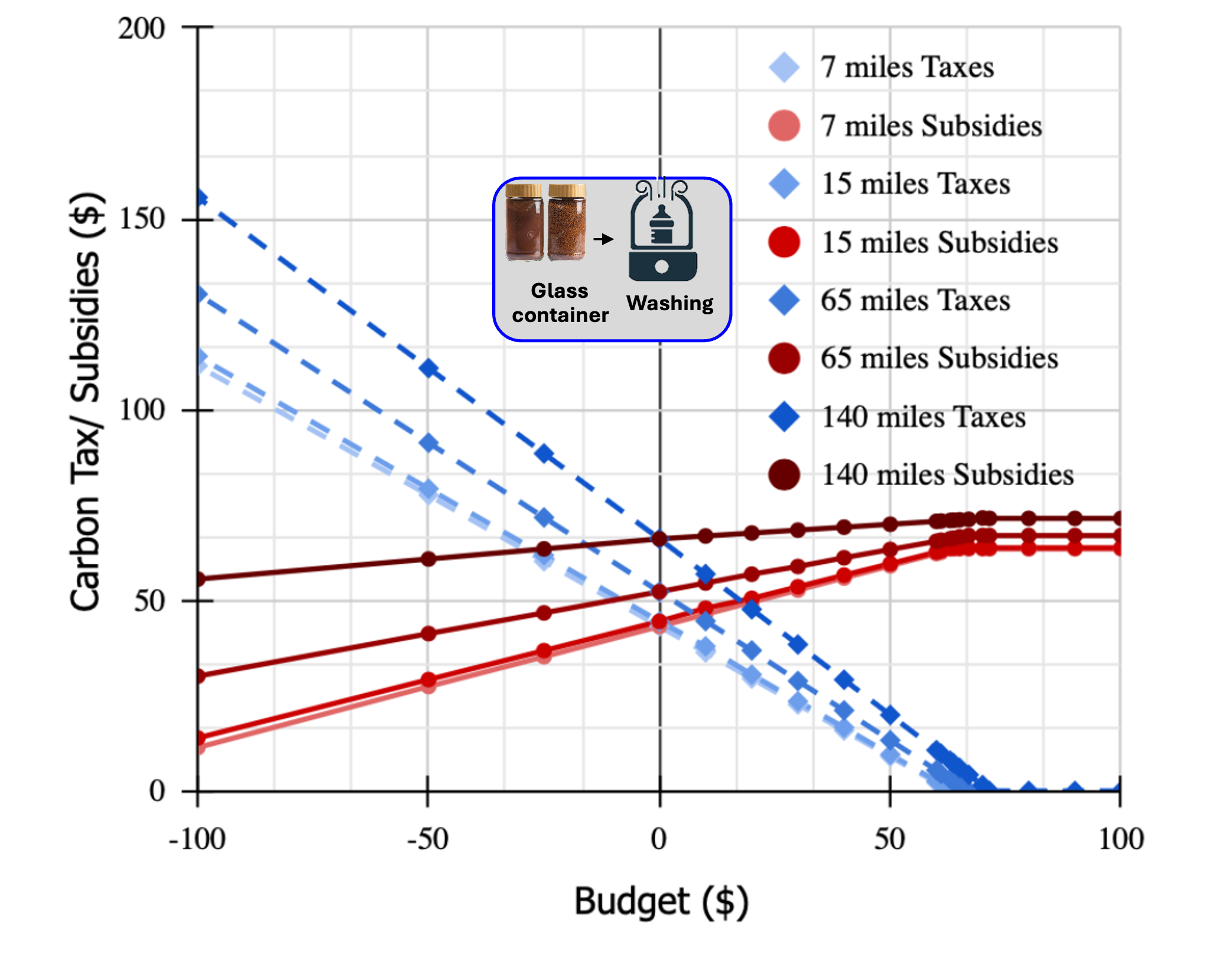}
    \caption{Sensitivity analysis of maximizing circularity for glass loss rates of 1\%, 3.13\%, and 10\%. (a) Carbon tax revenue and subsidies as functions of government budget. (b) Carbon tax rate versus government budget.}
    \label{fig:max Circularity gloss}
\end{figure}

When maximizing circularity, the scenario with 1\% to 10\% glass loss selects the \emph{glass} $\rightarrow$ \emph{glass washing} pathway. The relationship between the carbon tax, subsidies, and carbon tax rate with respect to the budget is illustrated in Fig.~\ref{fig:max Circularity gloss}. Similar to previous findings from the glass washing distance sensitivity analysis, an increase in glass loss percentage results in higher emissions, which in turn leads to a larger carbon tax and a steeper carbon tax slope. For the 1 and 3.13\% of glass loss, it can be seen that the subsidies provided increase while increasing the governmental budget. For a 10\% glass loss, the slope becomes so steep that the corresponding subsidy line exhibits a negative slope. This implies that an increase in the government budget would reduce the subsidies received by the industry, a result that appears counter-intuitive. However, when considering the overall cost of the industry, defined as the operating cost plus the carbon tax minus subsidies, although the subsidies decrease, the carbon tax decreases even more. Therefore, the overall cost still decreases as the government budget increases.

Furthermore, an increase in glass loss might intuitively suggest a shift toward alternative technologies. However, such a transition becomes both impractical and economically challenging due to the substantial subsidies required and the corresponding carbon taxes imposed. This is illustrated by the fact that the 10\% glass loss curve lies above the other solutions in Fig.~\ref{fig:max Circularity gloss}a, and the carbon tax rate lies above the other solutions in Fig.~\ref{fig:max Circularity gloss}b. Therefore, glass losses exceeding 10\% were not considered in this analysis.

The results obtained from the bilevel framework provide actionable insights for policy design by linking budget allocation to policy thresholds and technology adoption outcomes. Policymakers can use the proposed formulation to evaluate taxes and subsidies both independently and jointly across sectors, identifying the magnitude of policy instruments required to incentivize technologies that achieve higher circularity or greater emissions reductions. By explicitly accounting for budget constraints, the framework allows decision-makers to not only understand the policies required, but also assess the level of financial support needed to meet different sustainability objectives. Moreover, given the dynamic nature of the plastic waste management sector, the framework facilitates the evaluation of policy robustness under changing system parameters, supporting the design of policies that remain effective under evolving market and technological conditions.

\section{Conclusion}

This study presents a bilevel optimization framework that integrates carbon tax and subsidy policies to promote sustainability in supply chains, using a coffee packaging supply chain as a case study. Different bilevel formulations are considered based on two different environmental objectives: minimization of greenhouse gas (GHG) emissions and maximization of circularity, both subject to a governmental budget.

The results for single policy instruments demonstrate that subsidies are able to drive shifts toward the minimization of GHG emissions and the maximization of circularity, whereas a carbon tax alone can only induce shifts toward the minimization of emissions if imposed on a high enough level. This finding highlights that a carbon tax by itself struggles to transform the supply chain when the sustainability goal extends beyond GHG emissions alone. 

Combined policy instruments, including both carbon tax and subsidies, can effectively promote shifts toward technology pathways having the lowest GHG emissions and highest circularity. By solving the bilevel optimization problem, the relationship between governmental budget and the corresponding levels of carbon tax and subsidies can be determined. Assuming a high benchmark carbon tax rate in 2025 (\$0.1/kg-CO$_{2,\mathrm{eq}}$) within the carbon tax rate-budget relationship profile, we observe that the carbon tax plays a minor role in shifting the solutions compared to subsidies. This suggests that the current carbon tax rates are insufficient to effectively promote sustainability within this low-carbon-intensity supply chain.

We further demonstrate how critical factors such as glass washing distance and glass loss rates influence the selection of policy instruments and alter the carbon tax–subsidies–budget relationship through sensitivity analysis. The bilevel optimization framework enables quantification of the required adjustments in carbon tax rates, subsidies, and budget levels to drive technological transitions under the changes of these critical factors.

A limitation of the current model is the assumption of the unbounded carbon tax rates and fixed technological parameters, which may limit the scope of policy dynamics captured. For instance, a mix of end-of-life technology pathways would likely be selected if a carbon tax cap were applied. Future research could incorporate such tax caps, extend the approach to other supply chains, and explore more complex multi-stakeholder frameworks, including considerations of uncertainty. These extensions will enhance the applicability of bilevel optimization in supporting transitions toward a circular economy.

\section*{Declaration of generative AI and AI-assisted technologies in the writing process.}
During the preparation of this work, the authors used Chat-GPT in order to improve the readability and language of the manuscript. After using this tool/service, the authors reviewed and edited the content as needed and took full responsibility for the content of the published article.

\section*{Acknowledgments}
This material is based upon work supported by the U.S. Department of Energy, Office of Energy Efficiency and Renewable Energy, Bioenergy Technologies Office under Award Number DEEE0009285, and the Department of Chemical and Biological Engineering at the University of Wisconsin-Madison.

\bibliography{references}

@article{avraamidou2020circular,
  title={Circular Economy-A challenge and an opportunity for Process Systems Engineering},
  author={Avraamidou, Styliani and Baratsas, Stefanos G and Tian, Yuhe and Pistikopoulos, Efstratios N},
  journal={Computers \& Chemical Engineering},
  volume={133},
  pages={106629},
  year={2020},
  publisher={Elsevier}
}

@article{avraamidou2017multiparametric,
  title={A multiparametric mixed-integer bi-level optimization strategy for supply chain planning under demand uncertainty},
  author={Avraamidou, Styliani and Pistikopoulos, Efstratios N},
  journal={IFAC-PapersOnLine},
  volume={50},
  number={1},
  pages={10178--10183},
  year={2017},
  publisher={Elsevier}
}

@techReport{Sagasta2017,
   author = {Javier Mateo-Sagasta and Sara Marjani Zadeh and Hugh Turral},
   institution = {FAO, IWMI},
   title = {Water pollution from agriculture: a global review Executive summary LED BY},
   year = {2017}
}

@article{Teixeira2025,
   author = {Antoine Teixeira and Julien Lefèvre},
   doi = {10.1111/jiec.70001},
   issn = {1088-1980},
   issue = {3},
   journal = {Journal of Industrial Ecology},
   month = {6},
   pages = {733-745},
   title = {Global supply chains and domestic climate policy: Addressing the substantial material‐related carbon footprint of final consumption in France},
   volume = {29},
   year = {2025}
}

@article{Jin2014,
   author = {Mingzhou Jin and Nelson A. Granda-Marulanda and Ian Down},
   doi = {10.1016/j.jclepro.2013.08.042},
   issn = {09596526},
   journal = {Journal of Cleaner Production},
   month = {12},
   pages = {453-461},
   publisher = {Elsevier Ltd},
   title = {The impact of carbon policies on supply chain design and logistics of a major retailer},
   volume = {85},
   year = {2014}
}

@article{mackay2019self,
  title={Self-tuning networks: Bilevel optimization of hyperparameters using structured best-response functions},
  author={MacKay, Matthew and Vicol, Paul and Lorraine, Jon and Duvenaud, David and Grosse, Roger},
  journal={arXiv preprint arXiv:1903.03088},
  year={2019}
}

@misc{ellenmacarthur2024,
  author       = {{Ellen MacArthur Foundation}},
  title        = {Circular Economy: Introduction},
  year         = {2024},
  howpublished = {\url{https://www.ellenmacarthurfoundation.org/topics/circular-economy-introduction/overview}},
  note         = {Accessed: 2024-12-26}
}

@inproceedings{milyani2018optimally,
  title={Optimally designed subsidies for achieving carbon emissions targets in electric power systems},
  author={Milyani, Ahmad H and Kirschen, Daniel S},
  booktitle={2018 IEEE Power \& Energy Society General Meeting (PESGM)},
  pages={1--5},
  year={2018},
  organization={IEEE}
}

@article{camacho2024metaheuristics,
  title={Metaheuristics for bilevel optimization: A comprehensive review},
  author={Camacho-Vallejo, Jos{\'e}-Fernando and Corpus, Carlos and Villegas, Juan G},
  journal={Computers \& Operations Research},
  volume={161},
  pages={106410},
  year={2024},
  publisher={Elsevier}
}

@article{stumpf2023circular,
  title={Circular plastics packaging--Prioritizing resources and capabilities along the supply chain},
  author={Stumpf, Lukas and Sch{\"o}ggl, Josef-Peter and Baumgartner, Rupert J},
  journal={Technological Forecasting and Social Change},
  volume={188},
  pages={122261},
  year={2023},
  publisher={Elsevier}
}

@article{Chen2018,
   author = {Wanting Chen and Zhi Hua Hu},
   doi = {10.1016/j.jclepro.2018.08.007},
   issn = {09596526},
   journal = {Journal of Cleaner Production},
   month = {11},
   pages = {123-141},
   publisher = {Elsevier Ltd},
   title = {Using evolutionary game theory to study governments and manufacturers’ behavioral strategies under various carbon taxes and subsidies},
   volume = {201},
   year = {2018}
}

@misc{Mikkonen2021Finland,
  author       = {Krista Mikkonen},
  title        = {How can Finland confront the climate crisis and preserve its biodiversity?},
  howpublished = {Open Access Government},
  month        = may,
  year         = {2021},
  url          = {https://www.openaccessgovernment.org/loss-of-biodiversity/111415/},
  note         = {Accessed: 2025‑07‑03}
}

@misc{EU_Directive2019_904,
  author       = {{European Parliament and Council of the European Union}},
  title        = {{Directive (EU) 2019/904 of 5 June 2019 on the reduction of the impact of certain plastic products on the environment}},
  howpublished = {Official Journal of the European Union (OJ L 155, pp. 1–19)},
  month        = jun,
  year         = {2019},
  note         = {Text with EEA relevance; accessed via EUR‑Lex},
  url          = {https://eur-lex.europa.eu/eli/dir/2019/904/oj},
  urldate      = {2025-07-07}
}

@article{CalistoFriant2021,
   author = {Martin Calisto Friant and Walter J.V. Vermeulen and Roberta Salomone},
   doi = {10.1016/j.spc.2020.11.001},
   issn = {23525509},
   journal = {Sustainable Production and Consumption},
   month = {7},
   pages = {337-353},
   title = {Analysing European Union circular economy policies: words versus actions},
   volume = {27},
   year = {2021}
}

@misc{EU_Commission2025_SUP,
  author       = {{European Commission}},
  title        = {Single‑Use Plastics},
  howpublished = {European Commission, Directorate‑General for Environment},
  year         = {2025},
  url          = {https://environment.ec.europa.eu/topics/plastics/single‑use‑plastics_en},
  urldate      = {2025‑07‑03},
}

@misc{ClarkeVivier2024_FrenchIncentives,
  author       = {Lara Clarke and Walfroy Vivier},
  title        = {New French Incentives for Green Buildings and Green Industry},
  howpublished = {JDSupra (Jones Day)},
  month        = mar,
  year         = {2024},
  url          = {https://www.jdsupra.com/legalnews/new-french-incentives-for-green-1681972/},
  note         = {Accessed: 2025-08-03}
}

@misc{Microsoft2025_RepairabilityIndexFrance,
  author       = {{Microsoft Corporation}},
  title        = {Repairability Index for France},
  howpublished = {Microsoft Surface Support},
  year         = {2025},
  url          = {https://support.microsoft.com/en-us/surface/repairability-index-for-france-8aa5a99c-b562-4260-811c-0589362ae161},
  note         = {Accessed: 2025-08-03}
}

@misc{HandbookGermany_WasteSeparation,
  author    = {{Handbook Germany: Together}},
  title     = {Waste separation in Germany},
  year      = {2025},
  url       = {https://handbookgermany.de/en/waste-separation},
  urldate   = {2025-08-03}
}

@misc{EEA2023_NetherlandsWastePrevention,
  author    = {{European Environment Agency}},
  title     = {Netherlands – Waste prevention country profile},
  year      = {2023},
  url       = {https://www.eea.europa.eu/themes/waste/waste-prevention/countries/2023-waste-prevention-country-fact-sheets/netherlands_waste_prevention_2023},
  urldate   = {2025-07-03}
}

@misc{ICAP_Germany_nETS2025,
  author       = {{International Carbon Action Partnership}},
  title        = {German national emissions trading system},
  year         = {2025},
  howpublished = {ICAP ETS Factsheet},
  url          = {https://icapcarbonaction.com/en/ets/german-national-emissions-trading-system},
  urldate      = {2025-06-03}
}

@misc{China_14thFiveYearPlan_2021,
  title        = {Outline of the People’s Republic of China 14th Five-Year Plan for National Economic and Social Development and Long-Range Objectives for 2035},
  author       = {Xinhua News Agency},
  translator   = {Etcetera Language Group, Inc.},
  editor       = {Ben Murphy (CSET Translation Lead)},
  howpublished = {\url{https://cset.georgetown.edu/wp-content/uploads/t0284_14th_Five_Year_Plan_EN.pdf}},
  note         = {English translation; original Chinese version passed by National People’s Congress, March 2021},
  year         = {2021},
}

@misc{Chile_MMA_EconomiaCircular_Plasticos,
  author       = {Ministerio del Medio Ambiente, Gobierno de Chile},
  title        = {Plásticos — Economía Circular},
  howpublished = {\url{https://economiacircular.mma.gob.cl/plasticos/}},
  note         = {Página oficial que describe la Ley 21.368 sobre plásticos de un solo uso y botellas plásticas; incluye fechas, obligaciones, fases de implementación},
  year         = {2024},
  url          = {https://economiacircular.mma.gob.cl/plasticos/},
  urldate      = {2025-09-28}
}

@misc{toronto_circular_food_innovators_fund,
  author       = {{City of Toronto}},
  title        = {Circular Food Innovators Fund},
  year         = 2024,
  url          = {https://www.toronto.ca/services-payments/grants-incentives-rebates/circular-food-innovators-fund/},
  note         = {Accessed: 2025-09-29}
}

@article{Tsou2012,
   author = {Yung-Shan Tsou and Hsiao-Fan Wang},
   doi = {10.1080/10170669.2012.684408},
   issn = {1017-0669},
   issue = {4},
   journal = {Journal of the Chinese Institute of Industrial Engineers},
   month = {6},
   pages = {226-236},
   title = {Subsidy and penalty strategy for a green industry sector by bi-level mixed integer nonlinear programming},
   volume = {29},
   year = {2012}
}

@article{ma2021infrastructures,
  title={Infrastructures for phosphorus recovery from livestock waste using cyanobacteria: Transportation, techno-economic, and policy implications},
  author={Ma, Jiaze and Tominac, Philip and Pfleger, Brian F and Zavala, Victor M},
  journal={ACS Sustainable Chemistry \& Engineering},
  volume={9},
  number={34},
  pages={11416--11426},
  year={2021},
  publisher={ACS Publications}
}

@article{torres2016design,
  title={Design of multi-actor distributed processing systems: A game-theoretical approach},
  author={Torres, Ana I and Stephanopoulos, George},
  journal={AIChE Journal},
  volume={62},
  number={9},
  pages={3369--3391},
  year={2016},
  publisher={Wiley Online Library}
}

@article{tominac2017game,
  title={A game theoretic framework for petroleum refinery strategic production planning},
  author={Tominac, Philip and Mahalec, Vladimir},
  journal={AIChE Journal},
  volume={63},
  number={7},
  pages={2751--2763},
  year={2017},
  publisher={Wiley Online Library}
}

@article{fortuny1981representation,
  title={A representation and economic interpretation of a two-level programming problem},
  author={Fortuny-Amat, Jos{\'e} and McCarl, Bruce},
  journal={Journal of the operational Research Society},
  volume={32},
  number={9},
  pages={783--792},
  year={1981},
  publisher={Taylor \& Francis}
}

@article{sinha2017review,
  title={A review on bilevel optimization: From classical to evolutionary approaches and applications},
  author={Sinha, Ankur and Malo, Pekka and Deb, Kalyanmoy},
  journal={IEEE transactions on evolutionary computation},
  volume={22},
  number={2},
  pages={276--295},
  year={2017},
  publisher={IEEE}
}

@article{hansen1992new,
  title={New branch-and-bound rules for linear bilevel programming},
  author={Hansen, Pierre and Jaumard, Brigitte and Savard, Gilles},
  journal={SIAM Journal on scientific and Statistical Computing},
  volume={13},
  number={5},
  pages={1194--1217},
  year={1992},
  publisher={SIAM}
}

@article{kleinert2021survey,
  title={A survey on mixed-integer programming techniques in bilevel optimization},
  author={Kleinert, Thomas and Labb{\'e}, Martine and Ljubi{\'c}, Ivana and Schmidt, Martin},
  journal={EURO Journal on Computational Optimization},
  volume={9},
  pages={100007},
  year={2021},
  publisher={Elsevier}
}

@article{saharidis2009resolution,
  title={Resolution method for mixed integer bi-level linear problems based on decomposition technique},
  author={Saharidis, Georges K and Ierapetritou, Marianthi G},
  journal={Journal of Global optimization},
  volume={44},
  pages={29--51},
  year={2009},
  publisher={Springer}
}

@article{koppe2010parametric,
  title={Parametric integer programming algorithm for bilevel mixed integer programs},
  author={K{\"o}ppe, Matthias and Queyranne, Maurice and Ryan, Christopher Thomas},
  journal={Journal of optimization theory and applications},
  volume={146},
  number={1},
  pages={137--150},
  year={2010},
  publisher={Springer}
}

@article{avraamidou2019multi,
  title={A multi-parametric optimization approach for bilevel mixed-integer linear and quadratic programming problems},
  author={Avraamidou, Styliani and Pistikopoulos, Efstratios N},
  journal={Computers \& Chemical Engineering},
  volume={125},
  pages={98--113},
  year={2019},
  publisher={Elsevier}
}

@inproceedings{li2006hierarchical,
  title={A hierarchical particle swarm optimization for solving bilevel programming problems},
  author={Li, Xiangyong and Tian, Peng and Min, Xiaoping},
  booktitle={International Conference on Artificial Intelligence and Soft Computing},
  pages={1169--1178},
  year={2006},
  organization={Springer}
}

@article{yue2017stackelberg,
  title={Stackelberg-game-based modeling and optimization for supply chain design and operations: A mixed integer bilevel programming framework},
  author={Yue, Dajun and You, Fengqi},
  journal={Computers \& Chemical Engineering},
  volume={102},
  pages={81--95},
  year={2017},
  publisher={Elsevier}
}

@ARTICLE{pymoo,
    author={J. {Blank} and K. {Deb}},
    journal={IEEE Access},
    title={pymoo: Multi-Objective Optimization in Python},
    year={2020},
    volume={8},
    number={},
    pages={89497-89509},
}

@article{tua2020reusing,
  title={Reusing glass bottles in Italy: A life cycle assessment evaluation},
  author={Tua, Camilla and Grosso, Mario and Rigamonti, Lucia},
  journal={Procedia CIRP},
  volume={90},
  pages={192--197},
  year={2020},
  publisher={Elsevier}
}

@article{munoz2025integrated,
  title={Integrated decision-making approach for the simultaneous design of food packaging and waste management technologies to achieve a Circular Economy},
  author={Munoz-Briones, Paola A and del Carmen Mungu{\'\i}a-L{\'o}pez, Aurora and Sanchez-Rivera, Kevin and Avraamidou, Styliani},
  journal={Computers \& Chemical Engineering},
  pages={109269},
  year={2025},
  publisher={Elsevier}
}

@book{von1952theory,
  title={The Theory of the Market Economy},
  author={von Stackelberg, H.},
  lccn={52004949},
  url={https://books.google.com/books?id=o3ceAAAAIAAJ},
  year={1952},
  publisher={William Hodge}
}

@article{Xu2025,
   author = {Zhuo Xu and Kevin Sanchez-Rivera and Charles Granger and Panzheng Zhou and Aurora del Carmen Munguia-Lopez and Ugochukwu M. Ikegwu and Styliani Avraamidou and Victor M. Zavala and Reid C. Van Lehn and Ezra Bar-Ziv and Steven De Meester and George W. Huber},
   doi = {10.1038/s44286-025-00247-1},
   issn = {2948-1198},
   issue = {7},
   journal = {Nature Chemical Engineering},
   month = {7},
   pages = {407-423},
   title = {Solvent-based plastic recycling technologies},
   volume = {2},
   year = {2025}
}

@article{li2023examining,
  title={Examining how government subsidies influence firms’ circular supply chain management: The role of eco-innovation and top management team},
  author={Li, Miaomiao and Cao, Guikun and Cui, Li and Liu, Xiaoquan and Dai, Jing},
  journal={International Journal of Production Economics},
  volume={261},
  pages={108893},
  year={2023},
  publisher={Elsevier}
}

@article{martelli2020optimization,
  title={Optimization of renewable energy subsidy and carbon tax for multi energy systems using bilevel programming},
  author={Martelli, Emanuele and Freschini, Marco and Zatti, Matteo},
  journal={Applied energy},
  volume={267},
  pages={115089},
  year={2020},
  publisher={Elsevier}
}

@article{chalmardi2019bi, title={A bi-level programming model for sustainable supply chain network design that considers incentives for using cleaner technologies}, author={Chalmardi, Mazyar Kaboli and Camacho-Vallejo, Jos{\'e}-Fernando}, journal={Journal of cleaner production}, volume={213}, pages={1035--1050}, year={2019}, publisher={Elsevier} }

@article{camacho2023hierarchized, title={A hierarchized green supply chain with customer selection, routing, and nearshoring}, author={Camacho-Vallejo, Jos{\'e}-Fernando and D{\'a}vila, D{\'a}maris and Nucamendi-Guill{\'e}n, Samuel}, journal={Computers \& Industrial Engineering}, volume={178}, pages={109151}, year={2023}, publisher={Elsevier} }

@inproceedings{ma2018incentive,
  title={An Incentive-based Bi-level optimization Model for Collaborative Green Product Line Design},
  author={Ma, Shuang and Chen, Songlin and Cai, Xiaotian},
  booktitle={2018 IEEE International Conference on Industrial Engineering and Engineering Management (IEEM)},
  pages={981--985},
  year={2018},
  organization={IEEE}
}

@article{he2023hierarchical,
  title={Hierarchical optimization of policy and design for standalone hybrid power systems considering lifecycle carbon reduction subsidy},
  author={He, Yi and Guo, Su and Dong, Peixin and Huang, Jing and Zhou, Jianxu},
  journal={Energy},
  volume={262},
  pages={125454},
  year={2023},
  publisher={Elsevier}
}

@article{luo2022evaluating,
  title={Evaluating the impact of carbon tax policy on manufacturing and remanufacturing decisions in a closed-loop supply chain},
  author={Luo, Ruiling and Zhou, Li and Song, Yang and Fan, Tijun},
  journal={International Journal of Production Economics},
  volume={245},
  pages={108408},
  year={2022},
  publisher={Elsevier}
}

@article{mesrzade2023bilevel,
  title={A Bilevel Model for Carbon Pricing in a Green Supply Chain Considering Price and Carbon-Sensitive Demand},
  author={Mesrzade, Pegah and Dehghanian, Farzad and Ghiami, Yousef},
  journal={Sustainability},
  volume={15},
  number={24},
  pages={16563},
  year={2023},
  publisher={MDPI}
}

@article{hong2017optimizing,
  title={Optimizing an emission trading scheme for local governments: A Stackelberg game model and hybrid algorithm},
  author={Hong, Zhaofu and Chu, Chengbin and Zhang, Linda L and Yu, Yugang},
  journal={International Journal of Production Economics},
  volume={193},
  pages={172--182},
  year={2017},
  publisher={Elsevier}
}

@article{rahmani2024competitive,
  title={A competitive bilevel programming model for green, CLSCs in light of government incentives},
  author={Rahmani, Arsalan and Hosseini, Meysam and Sahami, Amir},
  journal={Journal of Mathematics},
  volume={2024},
  number={1},
  pages={4866890},
  year={2024},
  publisher={Wiley Online Library}
}

@article{liao2023optimal,
  title={Optimal decisions of closed-loop supply chain under government recycling subsidy and value co-creation},
  author={Liao, Bifeng and Shang, Yuanchao and Sun, Yangyi},
  journal={Frontiers in Energy Research},
  volume={11},
  pages={1308800},
  year={2023},
  publisher={Frontiers Media SA}
}

@article{caselli2024bilevel,
  title={Bilevel optimization with sustainability perspective: a survey on applications},
  author={Caselli, Giulia and Iori, Manuel and Ljubi{\'c}, Ivana},
  journal={arXiv preprint arXiv:2406.07184},
  year={2024}
}

@article{sudusinghe2022supply,
  title={Supply chain collaboration and sustainability performance in circular economy: A systematic literature review},
  author={Sudusinghe, Jayani Ishara and Seuring, Stefan},
  journal={International Journal of Production Economics},
  volume={245},
  pages={108402},
  year={2022},
  publisher={Elsevier}
}

@article{panza2025role,
  title={The role of carbon tax in the transition from a linear economy to a circular economy business model in manufacturing},
  author={Panza, Luigi and Peron, Mirco},
  journal={Journal of Cleaner Production},
  pages={144873},
  year={2025},
  publisher={Elsevier}
}

@inproceedings{kennedy1995particle,
  title={Particle swarm optimization},
  author={Kennedy, James and Eberhart, Russell},
  booktitle={Proceedings of ICNN'95-international conference on neural networks},
  volume={4},
  pages={1942--1948},
  year={1995},
  organization={ieee}
}

@article{beykal2020domino,
  title={Domino: Data-driven optimization of bi-level mixed-integer nonlinear problems},
  author={Beykal, Burcu and Avraamidou, Styliani and Pistikopoulos, Ioannis PE and Onel, Melis and Pistikopoulos, Efstratios N},
  journal={Journal of Global Optimization},
  volume={78},
  pages={1--36},
  year={2020},
  publisher={Springer}
}

@article{dempe2005discrete,
  title={Discrete bilevel programming: Application to a natural gas cash-out problem},
  author={Dempe, Stephan and Kalashnikov, Vyacheslav and R{\i}os-Mercado, Roger Z},
  journal={European Journal of Operational Research},
  volume={166},
  number={2},
  pages={469--488},
  year={2005},
  publisher={Elsevier}
}

@article{fischetti2018use,
  title={On the use of intersection cuts for bilevel optimization},
  author={Fischetti, Matteo and Ljubi{\'c}, Ivana and Monaci, Michele and Sinnl, Markus},
  journal={Mathematical Programming},
  volume={172},
  number={1},
  pages={77--103},
  year={2018},
  publisher={Springer}
}

@article{Baratsas2022,
   author = {Stefanos G. Baratsas and Efstratios N. Pistikopoulos and Styliani Avraamidou},
   doi = {10.1016/j.compchemeng.2022.107697},
   issn = {00981354},
   journal = {Computers and Chemical Engineering},
   month = {4},
   publisher = {Elsevier Ltd},
   title = {A quantitative and holistic circular economy assessment framework at the micro level},
   volume = {160},
   year = {2022}
}
\end{document}